\documentclass[letterpaper,11pt]{amsart}

\usepackage{amsmath,amssymb,amsxtra,amsthm,hyperref}
\usepackage{tikz}



\newtheorem*{thm:main}{Theorem \ref*{thm.main}}
\newtheorem*{thm:main1}{Theorem \ref*{thm.main.all}}

\newcommand{\bh}[1] {\mathcal{B}(\mathcal{#1})}
\newcommand{\pref}[1] {(\ref{#1})}

\newcommand{\lspan} {\operatorname{span}}

\newtheorem{theorem}{Theorem}[section]
\newtheorem{corollary}[theorem]{Corollary}
\newtheorem{proposition}[theorem]{Proposition}    
\newtheorem{lemma}[theorem]{Lemma}

\theoremstyle{definition}
\newtheorem{definition}[theorem]{Definition}

\newtheorem{example}[theorem]{Example}
\newtheorem{remark}[theorem]{Remark}

\numberwithin{equation}{section}

\begin{document}

\title{Regular Dilation on Graph Products of $\mathbb{N}$}

\author{Boyu Li}
\address{Pure Mathematics Department\\University of Waterloo\\Waterloo, ON\\Canada \ N2L--3G1}
\email{b32li@math.uwaterloo.ca}
\date{\today}


\begin{abstract} We extended the definition of regular dilation to graph products of $\mathbb{N}$, which is an important class of quasi-lattice ordered semigroups. Two important results in dilation theory are unified under our result: namely, Brehmer's regular dilation on $\mathbb{N}^k$ and Frazho-Bunce-Popescu's dilation of row contractions. We further show that a representation of a graph product has an isometric Nica-covariant dilation if and only if it is $\ast$-regular. A special case of our result was considered by Popescu, and we studied the connection with Popescu's work.
\end{abstract}

\maketitle

\section{Introduction}

Since the celebrated Sz.Nagy dilation theorem that showed a contraction has an isometric dilation, people have been trying to generalize this beautiful result to several variables. Ando showed that a pair of commuting contractions has a commuting isometric dilation. However, a well-known example due to Parrott shows a triple of commuting contractions may fail to have a commuting isometric dilation. Many studies consider what kind of extra conditions we need to guarantee a Sz.Nagy type dilation. There are two seemingly opposite directions this paper aims to unify. 

On one direction, Brehmer \cite{Brehmer1961} showed that if we put some extra conditions on the family of commuting contractions, then not only do they have an isometric dilation, but the isometric dilation actually satisfies a stronger condition known as regularity. Brehmer's result has been recently generalized to representations of any lattice ordered semigroup by the author \cite{BLi2014}. 

In another direction, we may consider non-commutative variables. Suppose $T_1,\cdots,T_n$ are contractions that are not necessarily commuting. It is observed by Frazho-Bunce-Popescu \cite{Frazho1982,Bunce1984,Popescu1989} that if $T_1,\cdots,T_n$ forms a row contraction in the sense that $\sum_{i=1}^n T_i T_i^*\leq I$, then there exists an isometric dilation for $T_i$ that is a row contraction. 

These two directions are seemingly unrelated. Indeed, one requires the commuting contractions to satisfy a stronger condition, whereas the other deals with non-commutative contractions. However, in this paper, I will show both results are special cases of $\ast$-regular dilation on graph products of $\mathbb{N}$. 

Graph products of $\mathbb{N}$ are considered as important examples of quasi-lattice ordered semigroups in \cite{CrispLaca2002}. Isometric Nica-covariant representations on a quasi-lattice ordered group have been intensively studied in the past decade. However, contractive Nica-covariant representations have only recently been defined in the lattice ordered group case (see \cite{Fuller2013, DFK2014}). Lattice ordered groups are quite restrictive compared to quasi-lattice ordered groups, and many interesting quasi-lattice ordered semigroups (e.g. the free semigroup $\mathbb{F}_n^+$) are not lattice ordered. This leads to the question for which type of representation on a quasi-lattice ordered group has a minimal isometric Nica-covariant dilation. We answer this question for this special class of quasi-lattice ordered semigroups of graph products of $\mathbb{N}$ and establish that having a minimal isometric Nica-covariant dilation is equivalent to being $\ast$-regular.

Popescu \cite{Popescu1999} showed that having a minimal isometric Nica-covariant dilation for a special class operators is equivalent to a property for which he calls the property (P). We extend Popescu's property (P) to the larger class of operators that correspond to representations of graph products of $\mathbb{N}$, and show that the property (P) holds whenever the representation is $\ast$-regular. 

\section{Background}

Throughout this paper, an operator $T$ is understood as a bounded linear operator on a complex Hilbert space $\mathcal{H}$. It is called \emph{a contraction} if its operator norm $\|T\|\leq 1$, and an \emph{isometry} if $T^* T=I$. If there is a larger Hilbert space $\mathcal{K}\supseteq\mathcal{H}$, we say that an operator $W\in\bh{K}$ is a \emph{dilation} of $T$ if $P_\mathcal{H} W^n \big|_\mathcal{H} = T^n$ for all $n\geq 1$. A familiar result due to Sarason \cite{Sarason1966} states that in such case, $\mathcal{K}$ decomposes as $\mathcal{H}_-\oplus\mathcal{H}\oplus\mathcal{H}_+$, and with respect to this decomposition, $$W=\begin{bmatrix} * & 0 & 0 \\ * & T & 0 \\ * & * & * \end{bmatrix}.$$

In particular, we say that $W$ is an \emph{extension} of $T$ if $\mathcal{H}$ is invariant. In other words, $\mathcal{H}_+=\{0\}$, and with respect to $\mathcal{K}=\mathcal{H}_-\oplus\mathcal{H}$, $$W=\begin{bmatrix} * & 0 \\ * & T \end{bmatrix}.$$

Dually, $W$ is called a \emph{co-extension} for $T$ if $\mathcal{H}^\perp$ is invariant for $W$. In other words, $\mathcal{H}_-=\{0\}$, and with respect to $\mathcal{K}=\mathcal{H}\oplus\mathcal{H}^+$, $$W=\begin{bmatrix} T & 0 \\ * & * \end{bmatrix}.$$

The celebrated Sz.Nagy dilation states that a contraction $T\in\bh{H}$ has an isometric co-extension. Moreover, the isometric co-extension $W$ can be chosen to be minimal in the sense that $$\mathcal{K}=\overline{span}\{W^k h:k\geq 0, h\in\mathcal{H}\}.$$ 

We call $W$ a \emph{minimal isometric dilation} for $T$. 

There are several attempts to generalize Sz.Nagy's result to the multivariate context. Ando \cite{Ando1963} proved that a pair of commuting contractions $T_1,T_2$ has commuting isometric dilation. However, the generalization to three commuting contractions fails as Parrott \cite{Parrott1970} gave a counter-example where three commuting contractions fails to have commuting isometric dilations. Given $n$ contractions $T_1,\cdots,T_n$, when can we find certain isometric dilation for them? 

There are two approaches to this question I would like to discuss in this paper. The first one requires some extra conditions of $T_i$. Brehmer \cite{Brehmer1961} first considered the question when $T_i$ has a stronger version of isometric dilation, now known as regular dilation. It has since been studied by many authors \cite{Halperin1962, SFBook, Gaspar1997}. It has been generalized to product systems \cite{Solel2008, Shalit2010}, and more recently, to any lattice ordered semigroup \cite{BLi2014}. For $n=(n_1,\cdots,n_k)\in\mathbb{Z}^k$, denote $T^n=\prod_{i=1}^k T_i^{n_i}$. Also denote $n^+=(n_1^+,\cdots,n_k^+)$ and $n^-=(n_1^-,\cdots,n_k^-)$, where $x^+=\max\{0,x\}$ and $x^-=\max\{0,-x\}$. It is clear that $n=n^+-n^-$. 

\begin{definition}\label{df.regular} An isometric dilation $(W_i)$ for $(T_i)$ is called \emph{regular} if it has an additional property that for any $n\in\mathbb{Z}^k$, $$T^{*n^-} T^{n^+}=P_\mathcal{H} W^{*n^-} W^{n^+} \big|_\mathcal{H}.$$

Dually, it is called \emph{$\ast$-regular} if for any $n\in\mathbb{Z}^k$, $$T^{n^+}T^{*n^-} =P_\mathcal{H} W^{*n^-} W^{n^+} \big|_\mathcal{H}.$$
\end{definition}  

For every subset $J\subset\{1,\cdots,n\}$, we denote $T_J=\prod_{j\in J} T_j$. Brehmer shows a necessary and sufficient condition for $T_i$ to have a regular (or $\ast$-regular) dilation:

\begin{theorem}[Brehmer \cite{Brehmer1961}]\label{thm.Brehmer} A commuting $n$-tuple of contractions $T_1,\cdots,T_n$ has a regular dilation if and only if for every $V\subset\{1,2,\cdots,k\}$, we have,

\begin{equation}\label{eq.Brehmer.regular}
\sum_{U\subset V} (-1)^{|U|} T_U^* T_U  \geq 0
\end{equation}

Here $|U|$ denotes the cardinality of $U$. 

Dually, $T_i$ has a $\ast$-regular dilation if and only if for every $V\subset\{1,2,\cdots,k\}$, we have,

\begin{equation}\label{eq.Brehmer}
\sum_{U\subset V} (-1)^{|U|} T_U T_U^* \geq 0
\end{equation}
\end{theorem} 

\begin{remark} Due to Theorem \ref{thm.Brehmer}, the family $(T_i)$ has a regular dilation if and only if $(T_i^*)$ has a $\ast$-regular dilation. To make the notation more consistent with the row contraction condition in the Frazho-Popescu-Bunce dilation, we shall mostly consider the $\ast$-regular dilation from now on.
\end{remark} 

Condition \pref{eq.Brehmer.regular} is a much stronger condition than the usual contractive condition that we require for an isometric dilation. For example, given two commuting contractions $T_1, T_2$, Ando's theorem always give an isometric dilation for this pair. However, to have a $\ast$-regular dilation, Brehmer's condition is equivalent of saying $I-T_1T_1^*-T_2T_2^* +T_1T_2T_2^* T_1^*\geq 0$. This can be false for many pairs of commuting contractions. 


The other approach to generalize Sz.Nagy's dilation is kind of the opposite to Brehmer's approach. Instead of adding more conditions to commutative $T_i$, we can replace the commutative condition by the row contractive condition, which states that $\sum_{i=1}^n T_i T_i^*\leq I$. Here, $T_i$ are no longer required to be commuting. It is now known as the Frazho-Bunce-Popescu dilation that if $T_1,\cdots,T_n$ is row contractive, then they can be dilated to row contractive isometries $V_1,\cdots,V_n$. 

Row contractive isometries (often called row isometries) is a well-studied subject in the field of operator algebra. For example, the $C^*$-algebra generated by row isometry is the well-known Cuntz-Toeplitz algebra. If we require further that $\sum_{i=1}^n V_i V_i^*=I$, its $C^*$-algebra is the renowned Cuntz-algebra. The $WOT$-closed non-self-adjoint algebra generated by a row isometry is a free semigroup algebra.

At first glance, Brehmer's result is seemingly unrelated to the Frazho-Bunce-Popescu dilation. Indeed, Brehmer's result is related to representations of the commutative semigroup $\mathbb{N}^k$ while the Frazho-Bunce-Popescu dilation is related to representations of the non-commutative free semigroup $\mathbb{F}_k^+$. Though regular dilation has been generalized to a larger class of semigroups called lattice ordered semigroups, the lattice order is a very restrictive requirement. For example, the free semigroup $\mathbb{F}_n^+$, and most quasi-lattice ordered semigroups, are not lattice ordered. 

This paper makes a first attempt to define regular and $\ast$-regular dilations on a larger class of semigroups, known as graph products of $\mathbb{N}$. It is also called the graph semigroup, or the right-angled Artin monoid (see \cite{CrispLaca2002, Eilers2016}). It is an important class of quasi-lattice ordered semigroups, recently studied in \cite{CrispLaca2002} for its Nica-covariant covariant representations. The main result of this paper generalizes both Brehmer's theorem and the Frazho-Bunce-Popescu dilation to this context.

Throughout this paper, we let $\Gamma$ denote a countable simple graph with vertex set $\Lambda$ and edge set $E(\Gamma)$. In other words, the vertex set $\Lambda$ is a countable set, and every edge $e\in E(\Gamma)$ corresponds to two distinct vertices $i,j\in\Lambda$. Every edge is undirected, and we say that $i$ \emph{is adjacent to} $j$ if there is an edge $e=(i,j)\in E(\Gamma)$.

The graph product of $\mathbb{N}$, $P_\Gamma=\Gamma_{i\in\Lambda} \mathbb{N}$, is defined to be the unital semigroup generated by $\{e_i\}_{i\in\Lambda}$, with additional rules that $e_i,e_j$ commutes whenever $i$ is adjacent to $j$.

A representation $T:P_\Gamma\to\bh{H}$ is uniquely determined by its value on the set of generators $T_i=T(e_i)$ that satisfies $T_iT_j=T_jT_i$ whenever $i$ is adjacent to $j$. We extend the definition of regularity (see Definition \ref{df.regular1}) to representations on $P_\Gamma$ and show that $T$ is $\ast$-regular if and only if for every finite subset $W\subseteq\Lambda$:
\begin{equation}\label{eq.main1}
\sum_{\substack{U\subseteq W \\ U\mbox{ is a clique}}} (-1)^{|U|} T_U T_U^*\geq 0.
\end{equation} 

Here, a set $U\subseteq\Lambda$ is called a \emph{clique} if vertices in $U$ are pairwise adjacent to one another. This unifies Brehmer's regular dilation on $\mathbb{N}^k$ and Frazho-Bunce-Popescu's dilation of row contractions. Indeed, when $\Gamma$ is the complete graph, every subset $U\subseteq W$ is a clique. In such case, Condition \pref{eq.main1} is the same as the $\ast$-regular Condition \pref{eq.Brehmer} in the Brehmer's result. When $\Gamma$ contains no edge, the only cliques in $\Gamma$ are the singletons. In such case, Condition \pref{eq.main1} is equivalent of saying $\sum_{i\in W} T_i T_i^*\leq I$ and thus a row contraction as in the Frazho-Bunce-Popescu's dilation. 

We further consider the question of when a representation of the graph product semigroup $P_\Gamma$ has a minimal isometric Nica-covariant dilation. According to \cite{Gaspar1997}, a pair of commuting contractions has isometric Nica-covariant dilation if and only if they are $\ast$-regular. Frazho-Bunce-Popescu's result also shows a row contraction can be dilated to isometries with pair-wise orthogonal range, which corresponds to an isometric Nica-covariant dilation on the free semigroup. We show that the $\ast$-regular condition is equivalent of having an isometric Nica-covariant dilation. 

To summarize, we establishs the following equivalence:

\begin{theorem}\label{thm.main.all} Let $T:P_\Gamma\to\bh{H}$ be a representation. Then the following are equivalent:
\begin{enumerate}
\item $T$ is $\ast$-regular,
\item $T$ has a minimal isometric Nica-covariant dilation,
\item $T$ satisfies Condition \pref{eq.main1} for every finite subset $W\subseteq\Lambda$.
\end{enumerate}
\end{theorem}

\section{Graph Products}

Fix a simple graph $\Gamma$ with a countable vertex set $\Lambda$. Recall that a graph product of $\mathbb{N}$  is a unital semigroup $P_\Gamma=\Gamma_{i\in\Lambda}\mathbb{N}$, generated by generators $\{e_i\}_{i\in\Lambda}$ where $e_i,e_j$ commute whenever $i,j$ are adjacent in $\Gamma$. We also call $P_\Gamma$ the graph semigroup or the right-angled Artin monoid. It is also closely related to the Cartier-Foata monoid \cite{Heaps} where $e_i, e_j$ commute whenever $i,j$ are not adjacent.

We can similarly define the graph product of $\mathbb{Z}$, $G_\Gamma=\Gamma_{i\in\Lambda} \mathbb{Z}$. It is defined to be the free product of $\mathbb{Z}$ modulo the rule that elements in the $i$-th and $j$-th copies of $\mathbb{Z}$ commute whenever $(i,j)$ is an edge of $\Gamma$. $G_\Gamma$ is a group, which is also called the graph group or the right-angled Artin group. $G_\Gamma$ together with $P_\Gamma$ is an important example of a quasi-lattice ordered group that is studied by Crisp and Laca \cite{CrispLaca2002}. 

\begin{example}\label{ex.graphprod}[Examples of Graph Products]
\begin{enumerate}
\item Consider the complete graph $\Gamma$ that contains every possible edge $(i,j)$ $i\neq j$. The graph product $\Gamma_{i\in\Lambda} \mathbb{N}$ is equal to the abelian semigroup $\mathbb{N}_k^+$, since any two generators $e_i,e_j$ commute.
\item Consider the graph $\Gamma$ that contains no edges. The graph product $P_\Gamma=\Gamma_{i\in\Lambda} \mathbb{N}$ is equal to the free product $\mathbb{F}_k^+$.
\item Consider the following graph product associated with the graph in Figure \ref{fg.1}.

\begin{figure}[h]
\begin{tikzpicture}[scale=0.75]

\draw [line width=1pt] (-1,-1) -- (-1,1);
\draw [line width=1pt] (-1,-1) -- (1,-1);
\draw [line width=1pt] (1,1) -- (1,-1);
\draw [line width=1pt] (1,1) -- (-1,1);
\draw [line width=1pt] (1,1) -- (-1,-1);

\node at (-1,1) {$\bullet$};
\node at (1,-1) {$\bullet$};
\node at (-1,-1) {$\bullet$};
\node at (1,1) {$\bullet$};

\node at (-1.35,1) {1};
\node at (1.35,1) {2};
\node at (-1.35,-1) {4};
\node at (1.35,-1) {3};
\end{tikzpicture}
\caption{A simple graph of 4 vertices}
\label{fg.1}
\end{figure}
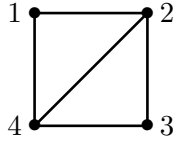

The graph product semigroup is a unital semigroup generated by $4$ generators $e_1,\cdots,e_4$, where the commutation relation is dictated by the edges of the graph. In this example, $e_i,e_j$ pairwise commute except for the pair $e_1, e_3$. 
\end{enumerate}
\end{example}

A typical element of $P_\Gamma$ is equal to $x=x_1x_2\cdots x_n$, where each $x_i$ belongs to a certain copy of $\mathbb{N}$. We often call $x$ an element of the semigroup, and $x_1x_2\cdots x_n$ an expression of $x$. Each $x_i$ is called \emph{a syllable} of this expression. We shall denote by $I(x_i)$ the index of $\mathbb{N}$ to which $x_i$ belongs. There might be many equivalent forms for $x$. First of all, $x$ may be able to be rewritten using fewer syllables: if $I(x_i)=I(x_{i+1})$, we may simply replace $x_ix_{i+1}=x_i'$, and write $x$ using $x_1\cdots x_{i-1} x_i' x_{i+2}\cdots x_n$. In particular, if $x_i=e$, we can always treat this $e$ as the identity for $I(x_{i+1})$ (or $I(x_{i-1})$), and thus $I(x_i)=I(x_{i+1})$. This process of merging two adjacent syllable from the same copy of $\mathbb{N}$ is called \emph{an amalgamation}. 

If $I(x_i)$ is adjacent to $I(x_{i+1})$ in the graph $\Gamma$, the commutative rules implies $x_ix_{i+1}=x_{i+1}x_i$, and thus $x$ can also be written as $$x_1\cdots x_{i-1}x_{i+1}x_i x_{i+2}\cdots x_n.$$  This process of switching two adjacent syllables when their corresponding copies of $\mathbb{N}$ commutes is called \emph{a shuffle}. We call two elements $x=x_1\cdots x_n$ and $y=y_1\cdots y_m$ \emph{shuffle equivalent} if we can obtain $y$ by repeatedly shuffling $x$. 

An expression $x=x_1\cdots x_n$ is called \emph{reduced} if it cannot be shuffled to another expression $x'$ which admits an amalgamation. Similar definitions of amalgamation and shuffle can be made for the graph group $\Gamma_{i\in\Lambda} \mathbb{Z}$. 

\begin{lemma}\label{lm.reduce} An expression $x=x_1\cdots x_n$ is reduced ($x$ can be in either $P_\Gamma$ or $G_\Gamma$) if and only if for all $i<j$ such that $I(x_i)=I(x_j)$, there exists an $i<t<j$ so that $I(x_t)$ is not adjacent to $I(x_i)$.
\end{lemma} 

The idea is that when $I(x_i)=I(x_j)$, as long as everything between $x_i$ and $x_j$ commute with $x_i$ and $x_j$, we can shuffle $x_j$ to be adjacent to $x_i$ and amalgamate the two. It is observed in \cite{Green1990} that reduced expressions are shuffle equivalent:

\begin{theorem}[Green \cite{Green1990}]\label{thm.shuffle} If $x=x_1\cdots x_n=x_1'\cdots x_m'$ are two reduced expressions for $x\in G_\Gamma$ (or $P_\Gamma$). Then two expressions are shuffle equivalent. In particular $m=n$.
\end{theorem}

This allows us to define \emph{the length of an element} $x$ to be $\ell(x)=n$, when $x$ has a reduced expression $x_1\cdots x_n$. 

Given a reduced expression $x=x_1\cdots x_n$, a syllable $x_i$ is called \emph{an initial syllable} if $x$ can be shuffled as $x=x_ix_1\cdots x_{i-1}x_{i+1}\cdots x_n$. Equivalently, it means the vertex $I(x_i)$ is adjacent to any previous vertices $I(x_j)$, $j<i$. The vertex $I(x_i)$ of an initial syllable is called \emph{an initial vertex}. The following lemma is partially taken from \cite[Lemma 2.3]{CrispLaca2002}. 

\begin{lemma}\label{lm.initial} Let $x=x_1\cdots x_n$ be a reduced expression. Then,
\begin{enumerate}
\item If $i\neq j$ and $x_i, x_j$ are two initial syllables, then $I(x_i)\neq I(x_j)$.
\item The initial vertices of $X$ are pairwise adjacent. 
\item Let $J=\{i: x_i\mbox{ is an initial syllable}\}$. Then $x=\prod_{j\in J} x_j \prod_{j\notin J} x_j$, where the second product is taken in the same order as in the original expression.
\end{enumerate}
\end{lemma} 

\begin{proof} If $I(x_i)=I(x_j)$ in a reduced expression, by Lemma \ref{lm.reduce}, there has to be an index $i<t<j$ so that $I(x_t)$ is not adjacent to $I(x_i)=I(x_j)$. Therefore, it is impossible to shuffle $x_j$ to the front. Therefore, any two initial syllables have different vertices. 

If $x_i,x_j$ are two initial syllables where $i<j$. Then to shuffle $x_j$ to the front, it must be the case that $x_j$ can commute with $x_i$, and thus $I(x_i)$ is adjacent with $I(x_j)$. This shows initial vertices are pairwise adjacent. 

Now let $J=\{1<j_1<j_2<\cdots<j_m\}$ be all $i$ where $x_i$ is an initial syllable. Then, we can recursively shift each $x_{j_s}$ to the front. The result is that we can shuffle all the initial vertices to the front as $\prod_{j\in J} x_j$, while all the other syllables are multiplied subsequently in the original order.
\end{proof} 

Lemma \ref{lm.initial} shows that the initial vertices are pairwise adjacent and thus form a clique of the graph $\Gamma$.

Lemma \ref{lm.initial} allows us to further divide a reduced expression of $x$ into blocks. Given a reduced expression $x=x_1\cdots x_n$, we define the first block $b_1$ of $x$ to be the product of all initial syllables. Since any two initial syllables commute, there is no ambiguity in the order of this product. We simply denote $I_1(x)=\{i: x_i\mbox{ is an initial syllable}\}$, and $b_1=\prod_{j\in I_1(x)} x_j$. Since $x_1$ is always an initial syllable, $I_1(x)\neq\emptyset$ and $b_1\neq e$. 

Now $x=b_1 x^{(1)}$, where $x^{(1)}$ has strictly shorter length compared to $x$. We can define the second block $b_2$ of $x$ to be the first block of $x^{(1)}$ when $x^{(1)}\neq e$. Of course, if $x^{(1)}=e$, we are finished since $x=b_1$. Repeat this process, and let each $x^{(t)}=b_{t+1} x^{(t+1)}$, where $b_{t+1}$ is the first block of $x^{(t)}$. Since the length of $x^{(t)}$ is always strictly decreasing, we eventually reach a state when $x^{(m-1)}=b_m x^{(m)}$ and $x^{(m)}=e$. In such case, $x$ is written as a product of $m$ blocks $x=b_1b_2\cdots b_{m}$. Here, each $b_j$ is the first block of $b_jb_{j+1}\cdots b_m$. We call this a block representation of $x$. We shall denote $I_t(x)$ be the vertex of all syllables in the $t$-th block $b_t$.

Since any two reduced expressions are shuffle equivalent, it is easy to see this block representation is unique. 

\begin{lemma}\label{lm.block} Let a reduced expression $x=x_1\cdots x_n$ have a block representation $b_1\cdots b_m$

\begin{enumerate}
\item\label{lm.block1} Two adjacent $I_t(x), I_{t+1}(x)$ are disjoint.
\item\label{lm.block2} For any vertex $\lambda_2\in I_{t+1}(x)$, there exists another vertex $\lambda_1\in I_t(x)$ so that $\lambda_1,\lambda_2$ are not adjacent.
\end{enumerate}
\end{lemma} 

\begin{proof} 
For \pref{lm.block1}, if $I_t(x),I_{t+1}(x)$ share some common vertex $\delta$, then the syllable corresponding to $\delta$ in the $(t+1)$-th block can be shuffled to the front of the $(t+1)$-th block, and since $\delta\in I_t(x)$, this syllable commutes with all syllable in the $t$-th block. Therefore, it can be amalgamated into the $t$-th block, leading to a contradiction that the expression is reduced.

For \pref{lm.block2}, if otherwise, we can pick a vertex $\lambda_2\in I_{t+1}(x)$ that is adjacent to every vertex in $I_t(x)$. The syllable corresponding to $\lambda_2$ can be shuffled to the front of $(t+1)$-th block, and commutes with everything in the $t$-th block. Therefore, it must be an initial syllable for $b_tb_{t+1}\cdots b_m$. But in such case, $\delta\in I_t(x)$ and cannot be in $I_{t+1}(x)$ by \pref{lm.block1}.
\end{proof} 

Studying regular dilations often requires a deep understanding of elements of the form $x^{-1}y$ for $x,y$ from the semigroup.

\begin{lemma}\label{lm.remove.ini} Let $x,y\in P_\Gamma$. Then, there exists $u,v\in P_\Gamma$ with $x^{-1}y=u^{-1}v$, and $I_1(u)$ is disjoint from $I_1(v)$. Moreover, $u,v$ are unique.
\end{lemma} 

\begin{proof} Suppose that there exists a vertex $\lambda\in I_1(x)\bigcap I_1(y)$. Then we can find initial syllables $e_\lambda^{m_1}$ and $e_\lambda^{m_2}$ from reduced expressions of $x,y$. We may without loss of generality assume that $x_1=e_\lambda^{m_1}$ and $y_1=e_\lambda^{m_2}$. 

Set $u_1=e_\lambda^{-\min\{m_1,m_2\}} x$ and $v_1=e_\lambda^{-\min\{m_1,m_2\}} y$. We have the relation $u_1^{-1} v_1=x^{-1}y$. Notice that at least one of $x_1$ and $y_1$ is removed in this process, and thus the total length $\ell(u_1)+\ell(v_1)$ is strictly less than $\ell(x)+\ell(y)$. Repeat this process whenever $I_1(u_j)\bigcap I_1(v_j)\neq \emptyset$, and recursively define $u_{j+1},v_{j+1}$ in the same manner to keep $u_j^{-1} v_j=u_{j+1}^{-1} v_{j+1}$. Since the total length $u_j,v_j$ is strictly decreasing in the process, we eventually stop in a state when $I_1(u_j)$ is disjoint from $I_1(v_j)$. This gives a desired $u=u_j, v=v_j$.

Suppose that $u^{-1}v=s^{-1} t$ for some other $s,t\in P_\Gamma$ with $I_1(s)\bigcap I_1(t)=\emptyset$. Let reduced expressions for $u,v,s,t$ be,
\begin{align*}
u &= u_1\cdots u_m \\
v &= v_1\cdots v_n \\
s &= s_1\cdots s_l \\
t &= t_1\cdots t_r 
\end{align*}

We first show $u^{-1}v=u_m^{-1}\cdots u_1^{-1} v_1 \cdots v_n$ is a reduced expression in $G_\Gamma$, and so is $s^{-1}t=s_l^{-1} \cdots s_1^{-1} t_1 \cdots t_r$. Assume otherwise, by Lemma \ref{lm.reduce}, there exists two syllables from the same vertex that commute with everything in between. These two syllables must have one from $u$ and the other from $v$, since $u_1\cdots u_m$ and $v_1\cdots v_n$ are both reduced. Let $u_i, v_j$ be two such syllables that come from the same vertex that commutes with everything in between. In that case, by Lemma \ref{lm.initial}, $u_i,v_j$ are both initial syllables for $u,v$. But $u,v$ have no common initial syllables, this leads to a contradiction. 

Therefore, $u_m^{-1}\cdots u_1^{-1} v_1 \cdots v_n=s_l^{-1} \cdots s_1^{-1} t_1 \cdots t_r$ are both reduced expressions for $u^{-1}v=s^{-1}t$, and thus by Theorem \ref{thm.shuffle} are shuffle equivalent. Notice each individual syllable $u_i, v_i, s_i, t_i$ is from the graph semigroup. To shuffle from $u_m^{-1}\cdots u_1^{-1} v_1 \cdots v_n$ to $s_l^{-1} \cdots s_1^{-1} t_1 \cdots t_r$, each $s_i^{-1}$ must be some $u_j^{-1}$, and $t_i$ must be some $v_j$. Therefore, $v_1\cdots v_n$ must be a shuffle of $t_1\cdots t_r$, and also $u_1\cdots u_m$ is a shuffle of $s_1\cdots s_l$. Hence, $s=u,t=v$. \end{proof} 

\begin{lemma}\label{lm.comm} Suppose $u,v\in\Gamma_{i\in\Lambda} \mathbb{N}$. Then the following are equivalent:
\begin{enumerate}
\item\label{lm.comm1} $u,v$ commute.
\item\label{lm.comm2} Every syllable $v_j$ of $v$ commutes with $u$.
\end{enumerate} 
\end{lemma} 

\begin{proof} 

\pref{lm.comm2}$\Longrightarrow$\pref{lm.comm1} is trivial. Assuming \pref{lm.comm1} and let $v=v_1\cdots v_m$. Consider the first syllable $v_1$ of $v$. Since $uv=vu$, $v_1$ is a initial syllable of $uv$. Therefore, $v_1$ commutes with $u$. By canceling $v_1$, one can observe that $v_2\cdots v_m$ also commutes with $u$, and recursively each $v_j$ commutes with $u$. \end{proof}

\begin{lemma}\label{lm.ini} Suppose $p\in P_\Gamma$, $\lambda\in\Lambda$ so that $\lambda\notin I_1(p)$ and $e_\lambda$ does not commute with $p$. Let $x,y\in P_\Gamma$ and apply the procedure in the Lemma \ref{lm.remove.ini} to repeatedly remove common initial vertex of $e_\lambda x$ and $py$ until $(e_\lambda x)^{-1} py=u^{-1}v$ with $I_1(u)\bigcap I_1(v)=\emptyset$. Then $u,v$ do not commute. 
\end{lemma} 

\begin{proof} Let $p=p_1\cdots p_n$ be a reduced expression of $p$. By Lemma \ref{lm.comm}, there exists a smallest $i$ so that $e_\lambda$ does not commute with $p_i$. We first observe that none of $p_1,\cdots,p_{i-1}$ come from the vertex $\lambda$.  Otherwise, if some $p_s$ comes from the vertex $\lambda$, it must commute with every $p_1,\cdots,p_{i-1}$ as $e_\lambda$ does. Therefore, $p_s$ is an initial syllable and $\lambda\in I_1(p)$, which contradicts to our assumption. 

Let $p_i$ be a syllable corresponding to vertex $\lambda'$, where $\lambda'$ is certainly not adjacent to $\lambda$.

Consider the procedure of removing a common initial vertex for $u_0=e_\lambda x$ and $v_0=py$. At each step, we removed a common initial vertex $\lambda_i$ for $u_i, v_i$ and obtained $u_{i+1}^{-1} v_{i+1}=u_i^{-1} v_i$, until we reach $u_m=u,v_m=v$ that shares no common initial vertex. It is clear that $\lambda\notin I_1(v_0)$ and $\lambda'\notin I_1(u_0)$. 

Observe that $\lambda_0\neq \lambda'$ since $\lambda\in I_1(e_\lambda x)$ and $\lambda'$ cannot be an initial vertex of $e_\lambda x$. Therefore, the syllable $p_i$ remains in $u_1$ after the first elimination step, while no syllable before $p_i$ belongs to the vertex $\lambda$. Hence, $\lambda\notin I_1(v_1)$ and $\lambda'\notin I_1(u_1)$. Inductively, $\lambda\notin I_1(v_j)$ and $\lambda'\notin I_1(u_j)$, and thus $e_\lambda$ is still an initial syllable of $u$ and $p_i$ is still a syllable of $v$. Therefore, $u,v$ do not commute. \end{proof} 

\section{Completely Positive Definite Kernels}

The problem of finding an isometric dilation turns out to be equivalent to showing that a certain kernel satisfies a so-called completely positive definite condition. Structures of completely positive definite kernels are studied in \cite{Popescu1996,Popescu1999b}, and we shall restate some of the results to our context.

Let $P$ be a unital semigroup sitting inside a group $G$ so that $P\bigcap P^{-1}=\{e\}$. For our purpose, the unital semigroup is taken to be a graph product $P_\Gamma=\Gamma_{i\in\Lambda} \mathbb{N}$, which lives naturally inside $G_\Gamma=\Gamma_{i\in\Lambda} \mathbb{Z}$. A \emph{unital Toeplitz kernel }on $P$ is a map $K:P\times P\to\bh{H}$ with the property that $K(e,e)=I$, $K(p,q)=K(q,p)^*$, and $K(ap,aq)=K(p,q)$ for all $a,p,q\in P$.

We call such a kernel \emph{completely positive definite }if for any $p_1,\cdots,p_n\in P$ and $h_1,\cdots,h_n\in\mathcal{H}$, we have $$\sum_{i,j=1}^n \left\langle K(p_i, p_j) h_j, h_i\right\rangle\geq 0.$$

Equivalently, this is saying that the $n\times n$ operator matrix $\left[K(p_i,p_j)\right]$, viewed as an operator on $\mathcal{H}^n$, is positive. Alternatively, each unital Toeplitz kernel $K$ corresponds to a map $\tilde{K}:P^{-1}P\to\bh{H}$, where $\tilde{K}(p^{-1} q)=K(p,q)$ and $\tilde{K}(x^{-1})=\tilde{K}(x)^*$. We shall abbreviate unital completely positive definite Toeplitz kernel as completely positive definite kernel. 

Existence of a completely positive definite kernel is closely related to the existence of an isometric dilation. A classical result known as Naimark dilation theorem \cite{Naimark1943} can be restated as the following theorem (\cite[Theorem 3.2]{Popescu1999b}):

\begin{theorem}\label{thm.Naimark} If $K$ is a completely positive definite kernel on a unital semigroup $P$, then there exists a Hilbert space $\mathcal{K}\supset\mathcal{H}$ and an isometric representation $V:P\to\bh{K}$ so that $$K(p,q)=P_\mathcal{H} V(p)^* V(q)\big|_\mathcal{H} \mbox{ for all }p,q\in P.$$

Moreover, $V$ can be taken as minimal in the sense that 
$$\overline{\lspan}\{V(p)h: p\in P,h\in\mathcal{H}\}=\mathcal{K},$$

The minimal isometric representation $V$ is unique up to unitary equivalence. 
\end{theorem}

Notice that in Theorem \ref{thm.Naimark}, if we set $p=e$, we get $K(e,q)=P_\mathcal{H} V(q)\big|_\mathcal{H}$. Assume now that $T:P\to\bh{H}$ is a contractive representation. If we can find a completely positive definite kernel $K$ so that $K(e,q)=T(q)$ for all $q\in P$, then Theorem \ref{thm.Naimark} gives us an isometric representation $V$ so that $T(q)=P_\mathcal{H} V(q)\big|_\mathcal{H}$. In other words, $V$ is an isometric dilation for $T$. Therefore, we reach the following conclusion:

\begin{corollary} Let $T:P\to\bh{H}$ be a contractive representation, for which there exists a completely positive definite kernel $K$ so that $K(e,q)=T(q)$. Then $T$ has an isometric dilation $V:P\to\bh{K}$, which can be taken as minimal in the sense that $$\overline{\lspan}\{V(p)h: p\in P,h\in\mathcal{H}\}=\mathcal{K}.$$

In particular, each $V(p)$ is a co-extension of $T(p)$. 
\end{corollary}

Such a kernel $K$ may not always exist. Indeed, if $P=\mathbb{N}^3$, let $T$ send three generators to the three commuting contractions as in the Parrott's example \cite{Parrott1970}. Such $T$ can never have an isometric dilation and thus there is no completely positive definite kernel $K$ so that $K(e,q)=T(q)$. Even when $T$ has an isometric dilation, $K$ may be extremely hard to define explicitly. 

Let us now turn our attention to contractive representations on a graph product $P_\Gamma=\Gamma_{i\in\Lambda} \mathbb{N}$. This semigroup is the free semigroup generated by $e_1,\cdots,e_n$ with additional rules that $e_i e_j=e_j e_i$ whenever $(i,j)\in E(\Gamma)$. Therefore, a representation $T$ of $P_\Gamma$ is uniquely determined by its values on generators $T_i=T(e_i)$, where they have to satisfy $T_i T_j=T_j T_i$ whenever $(i,j)\in E(\Gamma)$.

Let us fix a contractive representation $T:P_\Gamma\to\bh{H}$. We start by finding an appropriate completely positive definite kernel for $T$. Suppose an isometric regular dilation $V$ for $T$ exists, Theorem \ref{thm.Naimark} implies that $K(p,q)=P_\mathcal{H} V(p)^* V(q)\big|_\mathcal{H}$. Therefore, if $I_1(p)\bigcap I_1(q)\neq\emptyset$, $p=e_j p',q=e_j q'$ can both starts with a syllable in the same copy of $\mathbb{N}$. But since $V$ is isometric, $$V(p)^*V(q)=V(p')^*V(e_j)^*V(e_j)V(q')=V(p')^*V(q').$$

Therefore, it suffices to first consider the case that $I_1(p)\bigcap I_1(q)=\emptyset$. If otherwise, Lemma \ref{lm.remove.ini} gives $u,v$ so that $I_1(u)\bigcap I_1(v)=\emptyset$ and $u^{-1}v=p^{-1}q$. In such case, we shall define $K(p,q)=K(u,v)$. 

\begin{definition}\label{df.kernel} Given a contractive representation $T$ of the graph product $\Gamma_{i\in\Lambda} \mathbb{N}$, we define \emph{the Toeplitz kernel $K$ associated with $T$} using the following rules:
\begin{enumerate}
\item\label{df.proc1} $K(p,q)=T(q) T(p)^*$ whenever $I_1(p)\bigcap I_1(q)=\emptyset$ and $p,q$ commute.
\item\label{df.proc2} $K(p,q)=0$ whenever $I_1(p)\bigcap I_1(q)=\emptyset$ and $p,q$ do not commute.
\item\label{df.proc3} If there exists a vertex $i=I_1(p)\bigcap I_1(q)\neq\emptyset$. Let $p=e_i p'$, $q=e_i q'$. Define $K(p,q)=K(p',q')$.
\end{enumerate}  
\end{definition} 

\begin{remark} We may observe that since $I_1(e)=\emptyset$, and $e$ commutes with any $q$. $K(e,q)=T(q)$ by \pref{df.proc1}. Therefore, if $K$ is completely positive definite, the isometric Naimark dilation $V$ will be a dilation for $T$.
\end{remark}

\begin{remark} It follows from Lemma \ref{lm.remove.ini} that one can recursively remove common initial vertices from $p,q$ using \pref{df.proc3}, until we end up with unique $u,v$ with $u^{-1}v=p^{-1}q$ and $I_1(u)\bigcap I_1(v)=\emptyset$. Therefore, the Definition \ref{df.kernel} is well-defined for all pairs of $p,q$. 
\end{remark}

One can verify that the kernel $K$ is indeed a Toeplitz kernel. In fact, it satisfies a stronger property.

\begin{lemma}\label{lm.toeplitz} If $p,q,x,y\in P_\Gamma$ satisfies $p^{-1}q=x^{-1}y$, then $K(p,q)=K(x,y)$.
\end{lemma} 

\begin{proof} Repeatedly removing common initial vertices for the pairs $p,q$ and $x,y$ using the procedure in Lemma \ref{lm.remove.ini}, we end up with $p^{-1}q=u^{-1}v$, $x^{-1}y=s^{-1}t$, where $u,v$ has no common initial vertex; $s,t$ has no common initial vertex. Then, $K(p,q)=K(u,v)$ and $K(x,y)=K(s,t)$. By Lemma \ref{lm.remove.ini}, $u=s,t=v$. Therefore, $K(p,q)=K(x,y)$. 
\end{proof} 

\begin{definition}\label{df.regular1} We say that $T$ is \emph{$\ast$-regular} if the Toeplitz kernel $K$ associated with $T$ as defined in Definition \ref{df.kernel} is completely positive definite. A Naimark dilation $V$ for this kernel $K$ is called a $\ast$-regular dilation for $T$. Dually, we say that $T$ is \emph{regular} if $T^\ast$ is $\ast$-regular. Here, $T^*(e_i)=T(e_i)^*$.
\end{definition} 

\begin{remark} Our definition of regular dilation is slightly different from that of Brehmer's. When the graph semigroup is the abelian semigroup $\mathbb{N}^k$, Brehmer defined $T$ to be regular if a kernel $K^*$ is completely positive definite, where $K^*$ is the Toeplitz kernel by replacing Condition \pref{df.proc1} by $K^*(p,q)=T(p)^*T(q)$. In general, the kernel $K^*$ is different from the kernel we defined in Definition \ref{df.kernel}. However, it turns out when the semigroup is the abelian semigroup $\mathbb{N}^k$, our definition of regular dilation (Definition \ref{df.regular1}) coincides with Brehmer's definition (Definition \ref{df.regular}). 

However, on a general graph semigroup, when the kernel $K^*$ is completely positive definite is hard to characterize. For example, when the graph $\Gamma$ contains no edge and the graph semigroup corresponds to the free semigroup, the only chance that $p,q$ commute and $I_1(p)\bigcap I_1(q)=\emptyset$ is when at least one of $p,q$ is $e$. Therefore, in such case, $K^*=K$ and $K^*$ is completelely positive definite whenever $K$ is. 

Our definition of regular dilation implies there are isometric dilations for $T_i^*$ and thus co-isometric extensions for $T_i$. This coincides with the literature on the dilation of row contractions: for example, dilations for column contractions considered by Bunce \cite{Bunce1984} can be thought as regular dilation on the free semigroup $\mathbb{F}_+^k$.
\end{remark}

The $\ast$-regular representations are precisely those with a certain minimal Naimark dilation due to Theorem \ref{thm.Naimark}.

\begin{theorem} $T:P_\Gamma\to\bh{H}$ is $\ast$-regular if and only if it has a minimal isometric Naimark dilation $V:P_\Gamma\to\bh{K}$ so that for all $p,q\in P_\Gamma$, $K(p,q)=P_\mathcal{H} V(p)^* V(q)\big|_\mathcal{H}$. 
\end{theorem} 

\begin{remark} Given a representation $T:P_\Gamma\to\bh{H}$, there might be kernels different from the kernel we defined in Definition \ref{df.kernel} that are also completely positive definite. For example, it is pointed out in \cite{Opela2006} that when $\Gamma$ is acyclic, $T$ always has a unitary dilation. By restricting to $\mathcal{H}$, such a unitary dilation defines a completely positive definite kernel that is generally different from the kernel we defined. Popescu \cite{Popescu1999b} has also considered many ways to construct completely positive definite kernels on the free semigroup. 
\end{remark}

The goal of the next two sections is to provide a necessary condition for $\ast$-regularity of a contractive representation of a graph semigroup, which turns out to be also a sufficient condition. We draw our inspiration from two special cases where the graph is the complete graph and where the graph is the empty graph.

\begin{example}\label{ex.kernel.brehmer} In the case when $\Gamma$ is a complete graph on $k$ vertices. The graph semigroup $P_\Gamma$ is simply the abelian semigroup $\mathbb{N}^k$. It forms a lattice ordered semigroup. Each element in this semigroup can be written as a $k$ tuple $(a_1,\cdots,a_k)$. Since this semigroup is abelian, the set of initial vertex is precisely $\{i:a_i\neq 0\}$. 

Two elements $p=(p_i),q=(q_i)$ have disjoint initial vertex sets if and only if at least one of $p_i,q_i$ is zero for all $i$. In the terminology of the lattice order, this implies the greatest lower bound $p\wedge q=e$. As it is first defined in \cite{Brehmer1961}, a representation $T:\mathbb{N}^k\to\bh{H}$ is called $\ast$-regular if the kernel $K(p,q)$ is completely positive definite. 

Brehmer's result (Theorem \ref{thm.Brehmer}) shows that $K$ is completely positive definite if and only if for every subset $V\subseteq \{1,2,\cdots,k\}$, $$\sum_{U\subseteq V} (-1)^{|U|} T_U T_U^*\geq 0.$$

Here $|U|$ is the cardinality of $U$, and $T_U=\prod_{i\in U} T(e_i)$ with the convention that $T_\emptyset=I$.
\end{example}

\begin{example}\label{ex.kernel.popescu} In the case when $\Gamma$ is a graph on $k$ vertices with no edge. The graph semigroup $\Gamma_{i\in\Lambda} \mathbb{N}$ is simply the free semigroup $\mathbb{F}_k^+$. Fix a contractive representation $T:\mathbb{F}_k^+\to\bh{H}$, which is uniquely determined by its value on generators $T_i=T(e_i)$. The Toeplitz kernel associated with $T$ defined in Definition \ref{df.kernel} is the same as the kernel considered in \cite{Popescu1996, Popescu1999b}, where it is shown that $K$ is completely positive definite if and only if $T$ is row contractive in the sense that $$I-\sum_{i=1}^k T_i T_i^*\geq 0.$$

It turns out the minimal Naimark dilation for $K$ in this case is also a row contraction, and thus proves the Frazho-Bunce-Popescu dilation.
\end{example}

Inspired by both Example \ref{ex.kernel.brehmer} and \ref{ex.kernel.popescu}, our first main result unifies the Brehmer's dilation and the Frazho-Bunce-Popescu dilation. Recall that a set of vertices $U\subseteq\Lambda$ is called a clique if the subgraph induced on $U$ is a complete subgraph. 

\begin{theorem}\label{thm.main} Let $T$ be a contractive representation of a graph semigroup $P_\Gamma$. Then, $T$ is $\ast$-regular if for every $W\subseteq\Lambda$,

\begin{equation}\label{eq.main}
\sum_{\substack{U\subseteq W \\ U\mbox{ is a clique}}} (-1)^{|U|} T_U T_U^*\geq 0.
\end{equation} 
\end{theorem}

\begin{remark} Condition \pref{eq.main} coincides with the condition in both Example \ref{ex.kernel.brehmer} and \ref{ex.kernel.popescu}. Indeed, when $\Gamma$ is a complete graph, any $U\subseteq V$ is a clique. When $\Gamma$ contains no edge, the only cliques in $\Gamma$ are singletons $\{i\}$.
\end{remark}

\section{Technical Lemmas}

Since we are dealing with positive definiteness of operator matrices, the following lemma, taken from \cite[Lemma 14.13]{NestAlgebra}, is extremely useful.

\begin{lemma}\label{lm.Davidson} If an operator matrix $\begin{bmatrix} A & B^* \\ B & C\end{bmatrix}\in\mathcal{B}(\mathcal{H}_1\oplus\mathcal{H}_2)$ is positive, then there exists an operator $X:\mathcal{H}_1\to\mathcal{H}_2$ so that $B=XA^{1/2}$. Moreover, if $B$ has this form, then the operator matrix is positive if and only if $C\geq X X^*$. 
\end{lemma}


\begin{lemma}\label{lm.tech1} Let $X,L\in\bh{H}$ and $X\geq 0$. Define an $n\times n$ operator matrix $$A_n=
\begin{bmatrix} 
X & XL^* & XL^{*2} & \cdots & XL^{*(n-1)} \\
LX & X & XL^* & \cdots & XL^{*(n-2)} \\
L^2X & LX & X & \ddots & \vdots \\
\vdots & \ddots & \ddots & \ddots & XL^* \\
L^{n-1}X & L^{n-2}X & \cdots & LX & X
\end{bmatrix}.$$

If $LXL^*\leq X$, then every $A_n$ is positive. 

\end{lemma} 

\begin{proof} Assuming $LXL^*\leq X$, we shall inductively show each $A_n$ is positive. Since the case when $n=1$, $A_1=X\geq 0$ is given. Suppose $A_n\geq 0$, and rewrite $A_{n+1}$ as $$A_{n+1}=\left[
\begin{array}{ccccc|c}
 & & & & & XL^{*n} \\
 & & & & & XL^{*(n-1)} \\
 & & A_n & & & \vdots \\
 & & & & & \vdots \\
 & & & & & XL^* \\ \hline
L^n X & L^{n-1}X & \cdots & \cdots & LX & X \\
\end{array}\right].$$

Now notice that the row operator $[L^nX,\cdots,LX]=[0,\cdots,0,L] A_n$. Therefore, by Lemma \ref{lm.Davidson}, $A_{n+1}\geq 0$ if $$[0,\cdots,0,L] A_n \begin{bmatrix} 0 \\ \vdots \\ 0 \\ L^*\end{bmatrix} \leq X.$$

Expand the left hand side gives $LXL^*\leq X$. 
\end{proof} 

\begin{corollary}\label{cor.tech1} The matrix $A_n$ defined in Lemma \ref{lm.tech1} is positive if and only if $A_0=X\geq 0$ and $A_1\geq 0$. 
\end{corollary} 

\begin{proof} Indeed, $A_1=\begin{bmatrix} X & X^{1/2} X^{1/2}L^* \\ LX^{1/2} X^{1/2} & X\end{bmatrix}\geq 0$ if and only if $X\geq 0$ and $\left(LX^{1/2}\right)\left(X^{1/2}L\right)=LXL^*\leq X$ by Lemma \ref{lm.Davidson}. This is sufficient for every $A_n\geq 0$ by Lemma \ref{lm.tech1}.
\end{proof}

We now turn our attention to the contractive representation $T$ of a graph semigroup $P_\Gamma=\Gamma_{i\in\Lambda} \mathbb{N}$. Throughout this section, we fix such a representation $T$ and its associated Toeplitz kernel $K$ defined in Definition \ref{df.kernel}. For two finite subsets $F_1,F_2\subset P_\Gamma$, where $F_1=\{p_1,\cdots,p_m\}$ and $F_2=\{q_1,\cdots,q_n\}$, we denote $K[F_1,F_2]$ to be the $m\times n$ operator matrix, whose $(i,j)$-entry is equal to $K(p_i,q_j)$. When $F_1=F_2$, we simply write $K[F_1]=K[F_1,F_1]$. Recall $K$ is completely positive definite if and only if for all finite subsets $F\subseteq P_\Gamma$, $K[F]\geq 0$. If $F$ is a collection of elements that may contain duplicates, we may similarly define $K[F]$. It turns out duplicated elements will not affect the positivity of $K[F]$. 

\begin{lemma}\label{lm.rep} Let $F=\{p_1,p_1,p_2,\cdots,p_m\}$ and $F_1=\{p_1, p_2,\cdots,p_m\}$. Then $K[F]\geq 0$ if and only if $K[F_1]\geq 0$.
\end{lemma} 

\begin{proof} Denote $F_2=\{p_2,\cdots,p_m\}$. We have, 
$$K[F]=\left[\begin{array}{c|cc} 
I & I & K[p_1, F_2] \\ \hline
I & I & K[p_1,F_2] \\
K[F_2,p_1] & K[F_2,p_1] & K[F_2] \end{array} \right].$$
Here, the lower right corner is $K[F_1]$.

By Lemma \ref{lm.Davidson}, $K[F]\geq 0$ if and only if $K[F_2,p_1] K[p_1 F_2]\leq K[F_2]$. By Lemma \ref{lm.Davidson} again, this happens if and only if $K[F_1]\geq 0$.
\end{proof} 

\begin{lemma}\label{lm.add.left} Let $F_1=\{p_1,\cdots,p_m\}$ and $F_2=\{q_1,\cdots,q_n\}$ and fix a vertex $\lambda\in\Lambda$ so that $\lambda$ is not an initial vertex for any of the $p_i$. Let $D(\lambda,F_1)$ be a diagonal $m\times m$ operator matrix whose $i$-th diagonal entry is equal to $T(e_\lambda)^m$ if $e_\lambda$ commutes with $p_i$ and $0$ otherwise. Then, $K[F_1,e_\lambda^m \cdot F_2]=D(\lambda,F_1)\cdot K[F_1,F_2]$.
\end{lemma} 

\begin{proof} This is essentially proving that $K(p_i, e_\lambda^m q_j)=T(e_\lambda)^m K(p_i,q_j)$ if $e_\lambda$ commutes with $p_i$ and $0$ otherwise. 

Assuming first that $e_\lambda$ commutes with $p_i$. Then $p_i^{-1} e_\lambda^m q_j = e_\lambda^m p_i^{-1} q_j$. A key observation here is that when this happens, $p_i$ contains no syllable from the vertex $\lambda$. Since $e_\lambda$ commutes with every syllable of $p_i$, if there is a syllable of $p_i$ from the vertex $\lambda$, it must be an initial syllable, which contradicts to our selection of $p_i$.

Repeatedly removing common initial vertices for $p_i,q_j$ using Lemma \ref{lm.remove.ini}, we end up with $p_i^{-1} q_j=u^{-1}v$, where $u,v$ have no common initial vertex. It follows from the Definition \ref{df.kernel} that $K(p_i,q_j)=K(u,v)$. Notice that $I_1(e_\lambda^m v)$ includes $\lambda$ and every vertex in $I_1(v)$ that is adjacent to $\lambda$. Moreover, we observed that $\lambda\notin I_1(u)$. Therefore, we have $I_1(e_\lambda^m v)\bigcap I_1(u)=\emptyset$.

Suppose $u,v$ commute. Then $p_i^{-1} e_\lambda^m p_j=e_\lambda^m v u^{-1}=u^{-1} e_\lambda v $. Therefore, by Lemma \ref{lm.toeplitz}, $K(p_i, e_\lambda^m q_j)=K(u, e_\lambda^m v)$.  Hence, in this case, 
$$K(u, e_\lambda v)=T(e_\lambda)^m T(v) T(u)^*=T(e_\lambda)^mK(u,v).$$ 

If $u,v$ does not commute, $e_\lambda^m v$ also does not commute with $u$. Therefore, $K(u,v)=K(u,e_\lambda v)=0$.

Assume now that $e_\lambda$ does not commute with $p_i$. Consider the procedure of removing common initial syllables in $p_i$ and $e_\lambda^m q_j$: since $\lambda$ is not an initial vertex of $p_i$, each step we have to cancel out a syllable from $p_i$ and $q_j$ that both commute with $e_\lambda^m$. After each step of removing a common initial vertex, we removed some syllable from $p_i$ that commute with $e_\lambda$. Since $\lambda$ is not an initial vertex of $p_i$, each step will not cancel out any $e_\lambda^m$. Eventually, we always end up with $p_i^{-1} q_j= u^{-1} e_\lambda^m v$, where $u,e_\lambda^m v$ do not share any common initial vertex.

By Lemma \ref{lm.comm}, some syllable in $p_i$ does not commute with $e_\lambda$. Since all the syllables that got canceled commute with $e_\lambda$, there has to be some syllable in the left over $u$ that does not commute with $e_\lambda$. Therefore, $u$ and  $e_\lambda^m v$ do not commute. Hence, $K(u,e_\lambda^m v)=0$. \end{proof} 

As an immediate corollary, 

\begin{corollary}\label{cor.tech.F} Let $F=\{p_1,\cdots,p_n\}$ be a finite subset of $P_\Gamma$, and $\lambda\in\Lambda$ is a vertex that is not an initial vertex for any of $p_i$. For every $m\geq 0$, denote $F_m=\bigcup_{j=0}^m e_\lambda^j \cdot F$. Then $K[F_m]\geq 0$ if and only if $K[F]\geq 0$ and $K[F_1]\geq 0$. 
\end{corollary}

\begin{proof} For each $i\leq j$, $K[e_\lambda^i F, e_\lambda^j F]=K[F,e_\lambda^{j-i} F]$. Let $D=D(\lambda,F)$ be the $n\times n$ diagonal operator matrix, whose $(i,i)$-entry is $T(e_\lambda)$ if $e_\lambda$ commutes with $p_i$ and $0$ otherwise. It follows from Lemma \ref{lm.add.left} that $K[F,e_\lambda^{j-i} F]=D^{j-i}K[F]$. Similarly, for each $i>j$, $$K[e_\lambda^i F, e_\lambda^j F]=K[e_\lambda^j F, e_\lambda^i F]^*=K[F]D^{*(i-j)}.$$ 

Therefore,
$$K[F_m]=
\begin{bmatrix} 
K[F] & K[F]D^* & K[F]D^{*2} & \cdots & K[F]D^{*m} \\
DK[F] & K[F] & K[F] D^* & \cdots & K[F]D^{*(m-1)} \\
D^2K[F] & DK[F] & K[F] & \ddots & \vdots \\
\vdots & \ddots & \ddots & \ddots & K[F]D^* \\
D^m K[F] & D^{m-1}X & \cdots & DK[F] & K[F]
\end{bmatrix}.$$

Corollary \ref{cor.tech1} can be applied so that $K[F_m]\geq 0$ if and only if $K[F]\geq 0$ and $K[F_1]\geq 0$.
\end{proof}

\begin{lemma}\label{lm.reduction} Let $F_1=\{p_1,\cdots,p_n\}$, $F_2=\{q_1,\cdots,q_m\}$ be finite subsets of $P_\Gamma$, and $\lambda\in\Lambda$ is a vertex that is not an initial vertex for any of $p_i$ nor $q_j$. Suppose that $e_\lambda$ commutes with every $q_j$, but not with any $p_i$. Denote,
\begin{align*}
F_0 &= F_1\bigcup F_2 \\
F &= e_\lambda\cdot\left(F_1\bigcup F_2\right)\bigcup\left(F_1\bigcup F_2\right) \\
 &= e_\lambda F_0\bigcup F_0 \\
F' &= e_\lambda\cdot F_2\bigcup F_1\bigcup F_2 
\end{align*}

Then, $K[F]\geq 0$ if and only if $K[F']\geq 0$.
\end{lemma} 

\begin{proof} Let $D$ denote an $m\times m$ diagonal operator matrix whose diagonal entries are all $T(e_\lambda)$. Repeatedly apply Lemma \ref{lm.add.left},
\begin{equation*}
K[F] = \begin{bmatrix} 
K[F_1] & K[F_1,F_2] & 0 & K[F_1,F_2] D^* \\
K[F_2,F_1] & K[F_2] & 0 & K[F_2] D^* \\
0 & 0 & K[F_1] & K[F_1,F_2] \\
D K[F_2,F_1] & DK[F_2] & K[F_2,F_1] & K[F_2]
\end{bmatrix}.
\end{equation*}

Denote the upper left $2\times 2$ corner by $X= \begin{bmatrix} K[F_1] & K[F_1,F_2] \\ K[F_2,F_1] & K[F_2]\end{bmatrix}$. It is clear that $X=K[F_0]$. Let $L$ be a $(n+m)\times (n+m)$ diagonal operator matrix, whose first $n$ diagonal entries are $0$, and the rest $m$ diagonal entries be $T(e_\lambda)$. Then, the lower left $2\times 2$ corner can be written as $LX$, and $K[F]=\begin{bmatrix} X & XL^* \\ LX & X\end{bmatrix}$.

Lemma \ref{lm.tech1} states that $K[F]\geq 0$ if and only if $X=K[F_0]\geq 0$ and $LXL^*\leq X$. Explicitly writing out $X-LXL^*$, we get,
\begin{equation}\label{eq.reduction1}
X-LXL^* = \begin{bmatrix} K[F_1] & K[F_1,F_2] \\ K[F_2,F_1] & K[F_2] - DK[F_2]D^* \end{bmatrix}.
\end{equation} 

Now consider $K[F']$:
\begin{equation}
K[F'] = \begin{bmatrix} 
K[F_2] & 0 & K[F_2] D^* \\
0 & K[F_1] & K[F_1,F_2] \\
DK[F_2] & K[F_2,F_1] & K[F_2]
\end{bmatrix}.
\end{equation}

Notice here $\begin{bmatrix} 0 \\ DK[F_2] \end{bmatrix} = \begin{bmatrix} 0 \\ D\end{bmatrix} K[F_2]$. By Lemma \ref{lm.Davidson}, $K[F']\geq 0$ if and only if $K[F_2]\geq 0$ and $$\begin{bmatrix} 0 \\ D\end{bmatrix} K[F_2] \begin{bmatrix} 0 & D^*\end{bmatrix} = \begin{bmatrix} 0 & 0 \\ 0 & D K[F_2] D^*\end{bmatrix} \leq \begin{bmatrix} K[F_1] & K[F_1,F_2] \\ K[F_2,F_1] & K[F_2]\end{bmatrix}.$$ 

This is precisely the condition required in Condition \pref{eq.reduction1}. Therefore, combing the results from above, $K[F]\geq 0$ if and only if $K[F']\geq 0$, $K[F_0]\geq 0$ and $K[F_2]\geq 0$. But notice $F_0,F_2$ are subset of $F'$, the later condition is equivalent to $K[F']\geq 0$. \end{proof} 

\section{Proof of The Main Result}

We prove the first main result (Theorem \ref{thm.main}) in this section. The goal is to show that for every finite $F=\{p_1,\cdots,p_n\}\subset P_\Gamma$, $K[F]\geq 0$ where $K$ is the Toeplitz kernel associated with a contractive representation $T:P_\Gamma\to\bh{H}$ that satisfies Condition \pref{eq.main}. 

The plan to prove the main result Theorem \ref{thm.main} is divided into 2 steps. In the first step, we define an order on finite subsets of $P_\Gamma$, and show that for each $F\subset P_\Gamma$, $K[F]\geq 0$ follows from $K[F']\geq 0$ for some $F'<F$ under this order. This allows us to make an induction along finite subsets of $P_\Gamma$. 

The base case of the induction turns out to be the case when every element in $F$ has precisely one block. The second step is to show for all such $F$, $K[F]\geq 0$. Inspired by \cite[Section 6]{BLi2014}, we shall then use an argument to show such $K[F]$ can be decomposed as $RR^*$ for some operator matrix $R$ explicitly. 

For the first step, we show that as long as $F$ contains some element that has more than 1 block, one can find another finite subset $F'\subset P_\Gamma$ so that $K[F]\geq 0$ if $K[F']\geq 0$. The key is then to show that this process of finding $F'$ will terminate after finitely many steps. 

\begin{definition}\label{df.operation} For each $\lambda\in\Lambda$, and $p\in P_\Gamma$, define $d_\lambda(p)$ to be:
\begin{enumerate}
\item\label{rm.op1} If $p=e_\lambda^{n_1} p'\in F$ where $e_\lambda$ does not commute with $p'$, then $d_\lambda(p)=\{p'\}$.
\item\label{rm.op2} If $p=e_\lambda^{n_1} p'\in F$ where $e_\lambda$ commutes with $p'$, then $d_\lambda(p)=\{e_\lambda p', p'\}$.
\item\label{rm.op3} If $\lambda$ is not an initial vertex of $p$ and $e_\lambda$ does not commute with $p$, then $d_\lambda(p)=\{p\}$.
\item\label{rm.op4} If $\lambda$ is not an initial vertex of $p$ and $e_\lambda$ commutes with $p$,, then $d_\lambda(p)=\{e_\lambda p, p\}$.
\end{enumerate}

For any finite set $F\subseteq P_\Gamma$, denote $d_\lambda(F)=\bigcup_{p\in F} d_\lambda(p)$. 
\end{definition} 

\begin{lemma}\label{lm.main.reduction} Let $F=\{p_1,\cdots,p_n\}\subset P_\Gamma$ with some $p_i$ containing at least $2$ blocks. Pick a $\lambda$ that is an initial vertex for some $p_i$, but $e_\lambda$ does not commute with $p_i$. 

Then $K[F]\geq 0$ if $K[d_\lambda(F)]\geq 0$. 

\end{lemma} 

\begin{proof} Without loss of generality, assume $p_1$ has at least two blocks. First of all, by Lemma \ref{lm.block}, there exists an initial vertex $\lambda$ of $p_1$ that is not adjacent to some vertex $\lambda'$ in the second block of $p_1$. Therefore, $e_\lambda$ does not commute with $p_1$. We fix this vertex $\lambda$, and reorder $p_1,\cdots,p_n$ so that $\lambda$ is an initial vertex for $p_1,\cdots,p_m$ but not $p_{m+1},\cdots,p_n$. 

Write $p_i=e_\lambda^{n_i}p_i'$ for all $1\leq i\leq m$. Denote $F_0=\{p_1',\cdots,p_m',p_{m+1},\cdots,p_n\}$. None of elements in $F_0$ has $\lambda$ as an initial vertex. Let $N=\max\{n_i\}$ and denote $F_N=\bigcup_{j=0}^N e_\lambda^j\cdot F_0$. It is clear that $F\subseteq F_N$, and thus $K[F]\geq 0$ if $K[F_N]\geq 0$. By Corollary \ref{cor.tech.F}, $K[F_N]\geq 0$ if and only if $K[F_1]\geq 0$ where $F_1=\left(e_\lambda\cdot F_0\right)\bigcup F_0$. 

We may further split $F_0$ into two subsets $F_0=C\bigcup N$, where $C=\{f\in F: f\mbox{ commutes with }e_\lambda\}$ and $N=\{f\in F: f\mbox{ does not commute with }e_\lambda\}$. Now apply Lemma \ref{lm.reduction}, $K[F_1]\geq 0$ if and only if $K[\left(e_\lambda\cdot C\right)\bigcup F_0]\geq 0$. Denote $$F'=\left(e_\lambda\cdot C\right)\bigcup F_0=\left(e_\lambda\cdot C\right)\bigcup C \bigcup N.$$ This proves that $K[F']\geq 0$ implies $K[F]\geq 0$. 

To see $F'=d_\lambda(F)$: fix an element $p_i\in F$ and consider $4$ possibilities:
\begin{enumerate}
\item If $p_i=e_\lambda^{n_1} p_i'\in F$ where $e_\lambda$ does not commute with $p_i'$, then $d_\lambda(p_i)=\{p_i'\}$ is contained in $N\subseteq F_0\subseteq F'$;
\item If $p_i=e_\lambda^{n_1} p_i'\in F$ where $e_\lambda$ commutes with $p_i'$, then $p_i'$ is an element of $C$ and thus $d_\lambda(p_i)=\{e_\lambda p_i', p_i'\}$ is contained in $\left(e_\lambda\cdot C\right)\bigcup C\subseteq F'$;
\item If $\lambda$ is not an initial vertex of $p_i$ and $e_\lambda$ does not commute with $p_i$, then $p_i$ is in the set $N$ and $d_\lambda(p_i)=\{p_i\}$ is contained in $N\subseteq F'$;
\item If $\lambda$ is not an initial vertex of $p_i$ and $e_\lambda$ commutes with $p_i$, then $p_i$ is in the set $C$ and $d_\lambda(p_i)=\{e_\lambda p_i, p_i\}$ is contained in $\left(e_\lambda\cdot C\right)\bigcup C\subseteq F'$.
\end{enumerate}

One can now observe that $F'=d_\lambda(F)$. This finishes the proof.\end{proof}

\begin{remark}\label{rm.operation} One may observe that due to \pref{rm.op2} and \pref{rm.op4}, the set $F'$ might be a larger set compared to $F$. The idea here is we removed $e_\lambda$ where it does not commute with some later syllables, this should make syllables of each element in $F'$ more commutative with one another. Therefore repeating this process will end up with an $F'$ where every element has only one block. This motivates the Definition \ref{df.bvs}. 
\end{remark}

\begin{definition}\label{df.bvs} For each element $p\in P_\Gamma$ with $m$ blocks, we define \emph{the block-vertex sequence} of $p$ to be $m$ sets of vertices $B_1(p),\cdots,B_m(p)$, where $B_1(p)=\{\lambda\in I_1(p):e_\lambda\mbox{ does not commute with }p\}$, and $B_j(p)=I_j(p)$ for all $2\leq j\leq m$. In other words, the $j$-th set is equal to the vertex set of $j$-th block of $p$, except for the first block, where we only include any vertex that does not commutes with the rest of the blocks. We also define $B_0(p)=\{\lambda\in I_1(p): e_\lambda\mbox{ commutes with }p\}$, the set of all initial vertices that are adjacent to every vertex that appears in $p$.

Define the block-vertex length of $p$ be $c(p)=\sum_{j=1}^m \left|B_j(p)\right|$. 
\end{definition}

\begin{remark}\label{rm.bv} In the case that $p$ has only one block, then every syllable is initial and thus commuting. In such case, $B_1(p)=\emptyset$ and $c(p)=0$. This is the only case when $c(p)=0$.

Also observe that for $p=e_{\lambda_1}^{m_1}\cdots e_{\lambda_n}^{m_n}$, the power $m_i\geq 1$ does not affect the block-vertex sequence of $p$. The only thing that matters is what kind of vertex appears in each block. 

In a reduced expression of $p$, each syllable uniquely corresponds to some vertex in one of $B_0(p),\cdots,B_m(p)$. Therefore, the length $\ell(p)=\sum_{j=0}^m \left|B_j(p)\right|$. The quantity $c(p)=\ell(p)-|B_0(p)|$ counts the number of syllables that do not commute with the rest. 
\end{remark}

\begin{lemma}\label{lm.bv.basic} Let $p\in P_\Gamma$ and $\lambda\in\Lambda$.
\begin{enumerate}
\item\label{lm.bv.basic1} If $\lambda\in B_1(p)$, and $p=e_\lambda^n p'$. Then $c(p')<c(p)$.
\item\label{lm.bv.basic2} If $e_\lambda$ commutes with $p$, then the block vertex sequence of any element in $d_\lambda(p)$ is the same as that of $p$. Here, $d_\lambda(p)$ is defined as in the Definition \ref{df.operation}.
\item\label{lm.bv.basic3} If $e_\lambda$ does not commute with $p$ and $\lambda$ is not an initial vertex of $p$, then the block vertex sequence of any element in $d_\lambda(p)$ is the same as that of $p$.
\end{enumerate}
\end{lemma} 

\begin{proof} For \pref{lm.bv.basic1}, every vertex in $B_0(p)$ is still in $B_0(p')$. Since we removed the syllable $e_\lambda^n$, $\ell(p')\leq \ell(p)-1$, it is observed by Remark \ref{rm.bv} that $c(p')<c(p)$.

For \pref{lm.bv.basic2}, there are two cases: either $\lambda\in B_0(p)$ or not. In the first case, write $p=e_\lambda^n p'$ and $d_\lambda(p)=\{p,p'\}$. Since we only removed an initial vertex that commutes with the rest of the word, $p'$ has the same block-vertex sequence as $p$. In the later case when $\lambda\notin B_0(p)$, $d_\lambda(p)=\{p,e_\lambda p\}$. Since $e_\lambda$ commutes with $p$, $\lambda$ will be added to $B_0(e_\lambda p)$ and thus will not change the block-vertex sequence of $e_\lambda p$. In any case, the block vertex sequence of any element in $d_\lambda(p)$ is the same as that of $p$.

For \pref{lm.bv.basic3}, $d_\lambda(p)=\{p\}$, and it is clear.
\end{proof} 

\begin{lemma}\label{lm.bv} If $p_1,p_2$ have the same block-vertex sequence, then so does every element of $d_\lambda(p_1),d_\lambda(p_2)$. 
\end{lemma} 

\begin{proof} If $\lambda\in B_1(p_1)=B_1(p_2)$, write $p_i=e_\lambda^{n_i} p_i'$ and $d_\lambda(p_i)=\{p_i'\}$. Then $p_i'$ is $p_i$ with the syllable $e_\lambda^{n_i}$ removed, and since $p_1,p_2$ have the same block-vertex sequence, $p_1',p_2'$ must also have the same block-vertex sequence. In any other case, by Lemma \ref{lm.bv.basic}, every element in $d_\lambda(p_i)$ has the same block-vertex sequence as $p_i$. 
\end{proof} 

\begin{definition} Let $F\subset P_\Gamma$ be a finite set. Define $c(F)=\sum c(f)$, where the summation is over all $f\in F$, but multiple elements with the same block-vertex sequence are only summed once.
\end{definition} 

\begin{lemma}\label{lm.bv.reduce} $c(d_\lambda(F))<c(F)$. 
\end{lemma} 

\begin{proof} Without loss of generality, let $f_1,\cdots,f_t$ have distinct block-vertex sequences while $f_{t+1},\cdots,f_n$ have the same block vertex sequence as some $f_i$, $1\leq i\leq t$, where $f_1=p_1=e_\lambda^{n_1} p_1'$ and $e_\lambda$ not commuting with $p_1'$. Then $c(F)=\sum_{i=1}^t c(p_i)$.

Now, from Lemma \ref{lm.main.reduction}, $\lambda\in B_1(p_1)$. Therefore, $d_\lambda(f_1)=\{p_1'\}$, and $c(p_1')<c(f_1)$. Now apply Lemma \ref{lm.bv}, the block-vertex sequence of each $d_\lambda(f_{t+1}),\cdots,d_\lambda(f_n)$ is the same as that of some $d_\lambda(f_1),\cdots,d_\lambda(f_t)$. Moreover, by Lemma \ref{lm.bv.basic}, $c(d_\lambda(f_i))\leq c(f_i)$. Therefore, since $d_\lambda(F)=\bigcup_{i=1}^n d_\lambda(f_i)$, we have, $$c(d_\lambda(F))\leq\sum_{i=1}^t c(d_\lambda(f_i))<\sum_{i=1}^t c(f_i)=c(F). \qedhere$$ \end{proof} 

To summarize the first step towards the proof of the main theorem,

\begin{proposition}\label{prop.main.step1} For every finite subset $F\subset P_\Gamma$, there exists finite subset $\tilde{F}\subset P_\Gamma$, where every element in $\tilde{F}$ contains exactly one block, and $K[F]\geq 0$ if $K[\tilde{F}]\geq 0$. 
\end{proposition} 

\begin{proof} We start with $F=F_0$ and repeatedly apply Lemma \ref{lm.main.reduction} to obtain $F_1=d_\lambda(F)$, $F_2=d_\lambda(F_1),\cdots$. Lemma \ref{lm.main.reduction} proves that $K[F_n]\geq 0$ if $K[F_{n+1}]\geq 0$. Lemma \ref{lm.bv.reduce} shows that $c(F_n)$ is a strictly decreasing integral sequence, and thus must stop at some $F_N=\tilde{F}$. If $c(\tilde{F})\neq 0$, some elements in $\tilde{F}$ has at least 2 blocks and Lemma \ref{lm.main.reduction} can still be applied to obtain another set $\tilde{F}'=d_\lambda(\tilde{F})$ with $c(\tilde{F}')<c(\tilde{F})$. Therefore, the last $F_N=\tilde{F}$ must have $c(F_N)=0$, which is equivalent of saying every element in $\tilde{F}$ contains exactly one block. It is also clear that $K[F]\geq 0$ if $K[F_N]\geq 0$. 
\end{proof} 

Our second step shall prove that for every finite subset $F$ where every element has exactly one block, $K[F]\geq 0$. Since $F$ only containly finitely many syllables, we may consider only the case when $\Gamma$ is a finite graph. If an element has exactly one block, then every syllable commutes with all other syllables, and thus their vertices corresponds to a clique in $\Lambda$. For a clique $U$, denote $e_U=\prod_{\lambda\in J} e_\lambda$. Since $U$ is a clique, there is no ambiguity in the order of this product. One exception to the definition is that we shall consider the empty set as a clique as well, and denote $e_\emptyset = e$. When $\Gamma$ is a finite graph, there are only finitely many cliques. Denote $F_c=\{e_U:U\mbox{ is a clique}\}$. The first lemma shows that it suffices to prove $K[F_c]\geq 0$. 

\begin{lemma}\label{lm.kfc} If $K[F_c]\geq 0$, then for any finite subset $F$ of $P_\Gamma$ whose elements all have one block, $K[F]\geq 0$. 
\end{lemma} 

\begin{proof} Suppose $F=\{p_1,\cdots,p_n\}$ contains an element $e_\lambda^{n} p'$ with $n\geq 2$, then reorder $p_1,\cdots,p_n$ so that $\lambda$ is an initial vertex for $p_1,\cdots,p_m$ but not $p_{m+1},\cdots,p_n$. Let $p_i'$ be the $p_i$ with the syllable corresponding to $\lambda$ removed. Let $F_0=\{p_1',\cdots,p_m',p_{m+1},\cdots,p_n\}$ and let $C\subseteq F_0$ be all elements that commute with $e_\lambda$. Lemma \ref{lm.main.reduction} proves that $K[F_0]\geq 0$ if $K[F']=K[\left(e_\lambda\cdot C\right)\bigcup F_0]\geq 0$. Since elements in $F_0$ contain exactly one block, and elements in $C$ commute with $F_0$, we have every element in $F'$ contains exactly one block.

Moreover, each syllable corresponding to the vertex $\lambda$ is $e_\lambda$. Repeat this process until we reach $\tilde{F}$ where for all $\lambda$, all syllables corresponding to $\lambda$ are $e_\lambda$. In such case, every element has the form $e_U$ for some clique $U$. It is clear that $\tilde{F}\subset F_c$ and thus if $K[F_c]\geq 0$, then $K[\tilde{F}]\geq 0$ and thus $K[F]\geq 0$.
\end{proof} 

To show $K[F_c]\geq 0$, it suffices to show $K[F_c]$ can be decomposed as $R_c R_c^*$. Following the technique outlined in \cite[Section 6]{BLi2014}, we can explicitly find such $R_c$. Moreover, under a certain ordering, $R_c$ can be chosen to be a lower triangular matrix, and can thus be viewed as a Cholesky decomposition of $K[F_c]$. This will be done in Proposition \ref{prop.kfc}, where we shall see where the conditions in Condition \pref{eq.main} come from. 

From Condition \pref{eq.main}, denote
\begin{equation}\label{eq.main.Z}
Z_V=\sum_{\substack{U\subseteq V \\ U\mbox{ is a clique}}} (-1)^{|U|} T_U T_U^*\geq 0.
\end{equation}

Here, $V$ is any subset of the vertex set $\Lambda$, and $T_U=T(e_U)$. Assuming Condition \pref{eq.main} holds true for a contractive representation $T$, each $Z_V\geq 0$ and we can thus take its square root $Z_V^{1/2}\geq 0$. 

\begin{definition}\label{df.neighborhood} For a clique $V$, we define the neighborhood of $V$, denoted by $N_V$, to be $$N_V=\{\lambda\in\Lambda:\lambda\notin V,\mbox{ and }\lambda\mbox{ is adjacent to every vertex in }V\}.$$ In particular, we define $N_\emptyset=\Lambda$.
\end{definition} 

\begin{lemma}\label{lm.clique.id} Fix a clique $F$, then $$\sum_{\substack{F\subseteq W \\ W\mbox{ is a clique}}} T_{W\backslash F} Z_{N_W} T_{W\backslash F}^* = I.$$
\end{lemma}

\begin{proof} Replace $Z_{N_W}$ using Equation \pref{eq.main.Z},
\begin{align*}
& \sum_{\substack{W\supseteq F \\ W\mbox{ is a clique}}} T_{W\backslash F} Z_{N_W} T_{W\backslash F}^* \\
= & \sum_{\substack{W\supseteq F \\ W\mbox{ is a clique}}} T_{W\backslash F} \left(\sum_{\substack{U\subseteq N_W\\ U\mbox{ is a clique}}} (-1)^{|U|} T_U T_U^*\right) T_{W\backslash F}^* \\
= & \sum_{\substack{W\supseteq F \\ W\mbox{ is a clique}}} \left(\sum_{\substack{U\subseteq N_W\\ U\mbox{ is a clique}}} (-1)^{|U|} T_{(U\bigcup W)\backslash F} T_{(U\bigcup W)\backslash F} ^* \right)
\end{align*}

Suppose $U\subseteq N_W$ is a clique, then every vertex of $U$ is adjacent to every vertex in $W$, and vertices in $U$ are adjacent to one another. Therefore, $U\bigcup W$ is also a clique. The converse is true as well: if $U\bigcup W$ is a clique where $U\bigcap W=\emptyset$, then $U\subseteq N_W$ is a clique. Hence, we can rearrange the double summation so that we first sum over all possible cliques $V=U\bigcup W$, and then sum over all possible $U$. For a fixed clique $V=U\bigcup W$, the set $W=V\backslash U$ and the only requirement is that $F\subseteq W$. Therefore, we only sum those $U$ so that $U\subseteq V\backslash F$. Rewrite the double summation as:
$$\sum_{\substack{V=U\bigcup W \\ V\mbox{ is a clique}}}\left(\sum_{\substack{U\subseteq V\backslash F \\ U\mbox{ is a clique}}} (-1)^{|U|} T_{V\backslash F} T_{V\backslash F}^*\right).$$

For a fixed clique $V=U\bigcup W$ where $U\bigcap W=\emptyset$, consider the inner summation over all clique $U\subseteq V\backslash F$. $|U|$ can take any value between $0$ and $|V\backslash F|$. Moreover, for a fixed size $|U|=k$, there are precisely ${|V\backslash F| \choose |U|}$ possibilities for $U$ where $U\subseteq V\backslash F$ with size $k$. 

Therefore, the coefficient for $T_{V\backslash F}T_{V\backslash F}^*$ where $V$ is a clique containing $F$, is equal to 
$$\sum_{j=0}^{|V\backslash F|} {|V\backslash F| \choose j} (-1)^j$$

This summation is equal to $1$ if $V=F$ and $|V\backslash F|=0$. Otherwise, this is equal to $(1-1)^{|V\backslash F|}=0$. This proves the double summation is equal to $T_{F\backslash F}T_{F\backslash F}^*=I$.
\end{proof}

We are now ready to show $K[F_c]\geq 0$. $K[F_c]$ is a $|F_c|\times |F_c|$ operator matrix, whose rows and columns are indexed by cliques $U,V$. Its $(U,V)$-entry is equal to $K[e_U, e_V]$. Eliminating common initial vertices, $K[e_U,e_V]=K[e_{U\backslash V}, e_{V\backslash U}]$. Now $e_{U\backslash V}$ commutes with $e_{V\backslash U}$ if and only if all vertices in $U\backslash V$ are adjacent to all vertices in $V\backslash U$. In other words, $U\bigcup V$ is a clique. Therefore, we have,
\begin{equation}\label{eq.kernel.clique}
K[e_U,e_V]=\begin{cases}T_{V\backslash U} T_{U\backslash V}^*, \mbox{ if }U\bigcup V\mbox{ is a clique;} \\ 0, \mbox{ otherwise}.\end{cases}
\end{equation}

Let $R_c$ be a $|F_c|\times |F_c|$ operator matrix, where 
\begin{equation}\label{eq.kernel.R}
R_c[U,W]=\begin{cases} T_{W\backslash U} Z_{N_W}^{1/2}, \mbox{ if }U\subseteq W \\ 0, \mbox{ otherwise}.\end{cases}
\end{equation}

\begin{proposition}\label{prop.kfc} $K[F_c]=R_c\cdot R_c^*$. In particular, $K[F_c]\geq 0$.
\end{proposition} 

\begin{proof} The $(U,V)$-entry for $R_c\cdot R_c^*$ is equal to $\sum_W R_c[U,W]R_c[V,W]^*$. 

If $U\bigcup V$ is not a clique, we cannot find a clique $W$ that contains both $U$ and $V$. Therefore, for every clique $W$, we cannot have both $U,V$ contained in $W$. By Equation \pref{eq.kernel.R}, this implies at least one of $R_c[U,W],R_c[V,W]$ is $0$. Hence, the $(U,V)$-entry for $R_c\cdot R_c^*$ is $0$, which agrees with the $(U,V)$-entry of $K[F_c]$ by Equation \pref{eq.kernel.clique}.

If $U\bigcup V$ is a clique, then $R_c[U,W]R_c[V,W]^*$ may be non-zero only when $W$ is a clique containing both $U,V$. Therefore, in such case,
\begin{align*}
& \sum_W R_c[U,W]R_c[V,W]^* \\
=& \sum_{U\bigcup V\subseteq W} R_c[U,W]R_c[V,W]^* \\
=& \sum_{U\bigcup V\subseteq W} T_{W\backslash U} Z_{N_W} T_{W\backslash V}^*  \\
=& T_{V\backslash U} \left(\sum_{U\bigcup V\subseteq W} T_{W\backslash (U\bigcup V) } Z_{N_W} T_{W\backslash (U\bigcup V)}^*\right) T_{U\backslash V}^* 
\end{align*}

The summation in the middle is equal to $I$ by Lemma \ref{lm.clique.id}, in which $F$ is the fixed clique $U\bigcup V$. This proves that the $(U,V)$-entry for $R_c\cdot R_c^*$ is equal to $T_{V\backslash U} T_{U\backslash V}^*=K[e_U,e_V]$ in this case. 

Therefore, we conclude that $K[F_c]=R_c\cdot R_c^*$ and $K[F_c]\geq 0$. 
\end{proof} 

\begin{remark} We can regard $R_c$ as a Cholesky decomposition of $K[F_c]$ by rearranging $R_c$ as a lower triangular matrix. We first notice that whenever $U$ contains more elements than $W$, $R_c[U,W]=0$. Moreover, when $|U|=|W|$, $U\subseteq W$ is equivalent to $U=W$. Therefore, $R_c[U,W]=0$ whenever $|U|\leq |W|$ and $U\neq W$. Therefore, if we rearrange $F_c$ according to the size of cliques (larger cliques come first), $R_c$ becomes a lower triangular matrix.
\end{remark}

\begin{example} Let us consider the graph product of $\mathbb{N}$ associated with the graph in Figure \ref{fg.2}:
\begin{figure}[h]
\begin{tikzpicture}[scale=0.75]

\draw [line width=1pt] (0,1) -- (-1,0);

\node at (-1,0) {$\bullet$};
\node at (1,0) {$\bullet$};
\node at (0,1) {$\bullet$};

\node at (-1.35,0) {1};
\node at (1.35,0) {3};
\node at (0,1.35) {2};
\end{tikzpicture}
\caption{A Simple Graph on 3 Vertices}
\label{fg.2}
\end{figure}
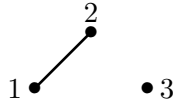

The graph semigroup is the unital semigroup generated by $e_1,e_2,e_3$ where $e_1,e_2$ commute. There are 5 cliques in this graph: $\{1,2\}$, $\{1\}$, $\{2\}$, $\{3\}$, and $\emptyset$. Under this ordering, 
$$K[F_c] = \begin{bmatrix} 
I & T_2^* & T_1^* & 0 & T_2^* T_1^* \\
T_2 & I & T_1^* T_2 & 0 & T_1^* \\
T_1 & T_2^* T_1 & I & 0 & T_2^* \\
0 & 0 & 0 & I & T_3^* \\
T_1 T_2 & T_1 & T_2 & T_3 & I
\end{bmatrix}.$$

We can write out the matrix $R_c$ using Equation \pref{eq.kernel.R}: 
$$R_c=\begin{bmatrix} 
I & 0 & 0 & 0 & 0 \\
T_2 & Z_{2}^{1/2} & 0 & 0 & 0 \\
T_1 & 0 & Z_1^{1/2} & 0 & 0 \\
0 & 0 & 0 & I & 0 \\
T_1 T_2 & T_1 Z_2^{1/2} & T_2 Z_1^{1/2} & T_3 & Z_{\{1,2,3\}}^{1/2}
\end{bmatrix}.$$ 

One can verify that $K[F_c]=R_c\cdot R_c^*$. \end{example}

We are now ready to prove the main Theorem \ref{thm.main}.

\begin{thm:main} Let $T$ be a contractive representation of a graph semigroup $P_\Gamma$. Then, $T$ is $\ast$-regular if for every finite $W\subseteq\Lambda$,

\begin{equation*}
\sum_{\substack{U\subseteq W \\ U\mbox{ is a clique}}} (-1)^{|U|} T_U T_U^*\geq 0.
\end{equation*} 
\end{thm:main}

\begin{proof} To show $T$ is $\ast$-regular given Condition \pref{eq.main}, it suffices to prove that the Toeplitz kernel $K$ in Definition \ref{df.kernel} is completely positive definite. For any finite subset $F\subset P_\Gamma$, it suffices to prove $K[F]\geq 0$. Proposition \ref{prop.main.step1} shows that it suffices to prove $K[\tilde{F}]\geq 0$ for some finite subset $\tilde{F}\subset P_\Gamma$, where each element in $\tilde{F}$ has precisely one block. Let $\Lambda_0$ be all the vertices that appears in a syllable of some element of $\tilde{F}$, which is a finite set. Denote $F_c=\{e_J\in\Lambda_0: J\mbox{ is a clique}\}$. By Lemma \ref{lm.kfc}, $K[\tilde{F}]\geq 0$ if $K[F_c]\geq 0$. Finally, by Proposition \ref{prop.kfc}, $K[F_c]\geq 0$. \end{proof} 

\begin{remark} The converse of Theorem \ref{thm.main} is also true (see Corollary \ref{cor.converse}).
\end{remark}

\section{Nica-Covariant Representation on Graph Products} 

Isometric Nica-covariant representations on a quasi-lattice ordered groups are first studied in \cite{Nica1992}, and were soon found to be an important concept in the study of operator algebras. Isometric Nica-covariant representations on graph semigroups, in particular graph products of $\mathbb{N}$, are intensively studied in \cite{CrispLaca2002}. It is observed in \cite[Theorem 24]{CrispLaca2002} that an isometric representation $V$ of the graph semigroup is isometric Nica-covariant if

\begin{enumerate}
\item for any two adjacent vertices $i,j$, $V_i$ and $V_j$ $\ast$-commute.
\item for any two non-adjacent vertices $i,j$, $V_i$ and $V_j$ have orthogonal ranges. In other words, $V_i^* V_j=0$. 
\end{enumerate}

Contractive Nica-covariant representations on lattice ordered semigroups are first defined and studied in \cite{Fuller2013, DFK2014}. However, lattice order is quite restrictive compared to quasi-lattice order. For example, the free semigroup $\mathbb{F}_m^+$ is quasi-lattice ordered, but not lattice ordered. In particular, the graph product $P_\Gamma$ is only lattice ordered when the graph $\Gamma$ is the complete graph, which corresponds to the abelian semigroup $\mathbb{N}^k$. This leads to a question of which representations of the graph product $P_\Gamma$ have isometric Nica-covariant dilations. 

In \cite{Gaspar1997}, it is shown that a pair of commuting contractions has a $\ast$-regular dilation if and only if they have a $\ast$-commuting isometric dilation, which is an equivalent way of saying a Nica-covariant dilation. The contractive Nica-covariant representations defined in \cite{Fuller2013, DFK2014, BLi2014} are always $\ast$-regular. It turns out that $\ast$-regular is equivalent of having an isometric Nica-covariant dilation.

\begin{theorem}\label{thm.nc} If $T:P_\Gamma\to\bh{H}$ is $\ast$-regular, then it has a minimal Naimark dilation that is an isometric Nica-covariant representation of the graph semigroup.
\end{theorem} 

The minimal Naimark dilation in Theorem \ref{thm.Naimark} can be constructed explicitly. We loosely follow the construction in \cite[Theorem 3.2]{Popescu1999b}. Given a completely positive definite kernel $K:P\times P\to\bh{H}$, define $\mathcal{K}_0=P\otimes\mathcal{H}$ with a semi-inner product defined by $$\left\langle \sum \delta_p\otimes h_p, \sum \delta_q\otimes k_q\right\rangle = \sum_{p,q} \langle K(q,p)h_p,k_q\rangle.$$

The original Hilbert space $\mathcal{H}$ can be embedded into $\mathcal{K}_0$ as $\delta_e\otimes\mathcal{H}$. The minimal Naimark dilation $V$ of $T$ acts on the $\mathcal{K}$ by $V(p)\delta_q\otimes h=\delta_{pq}\otimes h$, which are clearly isometries. Moreoever, for any $h_1,h_2\in\mathcal{H}$,
\begin{align*}
\langle V(q)^* V(p) h_1,h_2\rangle =& \langle \delta_p\otimes h_1, \delta_q\otimes h_2\rangle\\
=& \langle K(q,p) h_1, h_2\rangle
\end{align*}

Therefore, $P_\mathcal{H} V(q)^* V(p)\big|_\mathcal{H}=K(q,p)$. Let $\mathcal{N}=\{k\in\mathcal{K}_0: \langle k,k\rangle=0\}$. One can show that $\mathcal{N}$ is invariant for all $V(p)$, and thus we can let $\mathcal{K}=\overline{\mathcal{K}_0/\mathcal{N}}$, which is a Hilbert space. $V$ can be defined as isometries on $\mathcal{K}$, and it turns out that it is a minimal Naimark dilation. For technical details, one may refer to \cite[Theorem 3.2]{Popescu1999b}. It is worth noting that $\mathcal{H}$ is coinvariant for the minimal Naimark dilation $V$, and thus invariant for $V^*$.

To prove Theorem \ref{thm.nc}, it suffices to prove that the minimal Naimark dilation is isometric Nica-covariant. Throughout the rest of this section, we fix a contractive representation $T$ on $P_\Gamma$ that is $\ast$-regular, and let $V:P_\Gamma\to\bh{K}$ be the minimal Naimark dilation for $T$ described as above. 

\begin{lemma}\label{lm.orthogonal} Suppose $p\in P_\Gamma$, $\lambda\in\Lambda$ so that $\lambda\notin I_1(p)$ and $e_\lambda$ does not commute with $p$.  Then $V(e_\lambda)$ and $V(p)$ have orthogonal ranges. In other words, $V(e_\lambda)^* V(p)=0$. 
\end{lemma} 

\begin{proof} It suffices to prove for any $h=\sum_i \delta_{x_i}\otimes h_i\in\mathcal{K}_0=P_\Gamma\otimes\mathcal{H}$ and $h=\sum_j \delta_{y_j}\otimes k_i\in\mathcal{K}_0=P_\Gamma\otimes\mathcal{H}$, $\langle V(p) h, V(e_\lambda) k\rangle=0$. 

By the definition of the pre-inner product on $\mathcal{K}_0$, 
\begin{align*}
\langle V(p) h, V(e_\lambda) k\rangle &= \langle \sum_i \delta_{p\cdot x_i}\otimes h_i , \sum_j \delta_{e_\lambda \cdot y_j}\otimes k_i\rangle \\
&= \sum_{i,j} \langle K(e_\lambda \cdot y_j, p\cdot x_i) h_i, k_j\rangle
\end{align*}

Suppose $(e_\lambda \cdot y_j)^{-1} p\cdot x_i=u^{-1} v$ for some $u,v\in P_\Gamma$, where $u,v$ share no common initial vertices. By Lemma \ref{lm.ini}, $u,v$ do not commute. Therefore, $K(e_\lambda \cdot y_j, p\cdot x_i)=0$ for all $i,j$. Hence, the inner product is equal to $0$. \end{proof} 

\begin{lemma}\label{lm.nc} Let $p\in P_\Gamma$ and $\lambda\in\Lambda$ be a vertex such that $\lambda\notin I_1(p)$ and $e_{\lambda}$ commutes with $p$. Then $V(e_{\lambda})^* V(p)\big|_\mathcal{H}=V(p) V(e_{\lambda})^* \big|_\mathcal{H}$
\end{lemma} 

\begin{proof} By the minimality of $V$, $\lspan\{V(q)k:q\in P_\Gamma,k\in\mathcal{H}\}$ is dense in $\mathcal{K}$. Therefore, it suffices to prove for all $q\in P_\Gamma$, $h,k\in\mathcal{H}$,
\begin{equation}\label{eq.lm.nc}
\langle V(e_\lambda)^* V(p) h, V(q)k \rangle = \langle V(p) V(e_\lambda)^* h, V(q)k \rangle 
\end{equation}

Starting from the left hand side of Equation \pref{eq.lm.nc}, 
\begin{align*}
\langle V(e_\lambda)^* V(p) h, V(q)k \rangle =&  \langle V(e_\lambda q)^* V(p) h, k \rangle \\
=& \langle K(e_\lambda q,p) h, k \rangle \\
=& \langle K(q,p) T(e_\lambda)^* h, k \rangle
\end{align*}

Here we used Lemma \ref{lm.add.left} to show $K(e_\lambda q,p)=K(q,p) T(e_\lambda)^*$. Now since $V(e_\lambda)=\begin{bmatrix} T(e_\lambda) & 0 \\ * & * \end{bmatrix}$ with respect to the decomposition $\mathcal{K}=\mathcal{H}\oplus\mathcal{H}^\perp$, $V(e_\lambda)^* h = T(e_\lambda)^* h \in\mathcal{H}$. Therefore,
\begin{align*}
\langle K(q,p) T(e_\lambda)^* h, k \rangle =& \langle K(q,p) V(e_\lambda)^* h, k \rangle \\
=& \langle V(q)^* V(p) V(e_\lambda)^* h, k \rangle \\
=& \langle V(p) V(e_\lambda)^* h, V(q) k \rangle
\end{align*}

This proves Equation \pref{eq.lm.nc}. \end{proof} 

We now prove the main result of this section:

\begin{proof}[Proof of Theorem \ref{thm.nc}] It suffices to pick any two vertices $\lambda_1,\lambda_2$ and consider two cases when they are adjacent or not. 

If $\lambda_1,\lambda_2$ are not adjacent, by Lemma \ref{lm.orthogonal}, $V(e_{\lambda_1})$ and $V(e_{\lambda_2})$ are isometries with orthogonal ranges. 

If $\lambda_1,\lambda_2$ are adjacent, it suffices to prove for all $p\in P_\Gamma$, 
\begin{equation}\label{eq.nc1}
V(e_{\lambda_1})^* V(e_{\lambda_2}) V(p) \big|_\mathcal{H} = V(e_{\lambda_2}) V(e_{\lambda_1})^* V(p) \big|_\mathcal{H}. 
\end{equation}

Indeed, since $\lspan \{V(p)h:p\in P_\Gamma,h\in\mathcal{H}\}$ is dense in $\mathcal{K}$, Equation \pref{eq.nc1} implies that $V(e_{\lambda_1})^* V(e_{\lambda_2})= V(e_{\lambda_2}) V(e_{\lambda_1})^*$. 

There are now several possibilities:

If $\lambda\in I_1(p)$, we can write $p=e_{\lambda_1} p'$, and thus $V(p)=V(e_{\lambda_1})V(p')$. Since $\lambda_1,\lambda_2$ are adjacent, $V(e_{\lambda_1})$ commutes with $V(e_{\lambda_2})$. Hence, both sides of the Equation \pref{eq.nc1} are equal to $V(e_{\lambda_2}) V(p')\big|_\mathcal{H}$.

If $\lambda\notin I_1(p)$ and $e_{\lambda_1}$ does not commute with $p$, then $\lambda_1\notin I_1(e_{\lambda_2} p)$ and $e_{\lambda_1}$ does not commute with $e_{\lambda_2} p$ as well. Therefore, by Lemma \ref{lm.orthogonal}, $V(e_{\lambda_1})$ and $V(p)$ are isometries with orthogonal ranges, and $V(e_{\lambda_1})^* V(p)=0$. Similarly, $V(e_{\lambda_1})$ and $V(e_{\lambda_2} p)$ are isometries with orthogonal ranges, and $V(e_{\lambda_1})^* V(e_{\lambda_2} p)=0$. Both sides of the Equation \pref{eq.nc1} are $0$.

Lastly, if $\lambda\notin I_1(p)$ and $e_{\lambda_1}$ commutes with $p$. Then $e_{\lambda_2} p$ and $p$ are both element in $P_\Gamma$ that commutes with $e_{\lambda_1}$ without $\lambda_1$ as an initial vertex.  By Lemma \ref{lm.nc}, for every $h\in\mathcal{H}$, 
\begin{align*}
V(e_{\lambda_1})^* V(e_{\lambda_2}) V(p) h =&  V(e_{\lambda_2}) V(p) V(e_{\lambda_1})^* h\\
=& V(e_{\lambda_2}) V(e_{\lambda_1})^* V(p)  h
\end{align*} 

This is precisely the Equation \pref{eq.nc1}, and thus we finished the proof. \end{proof} 

\begin{corollary}\label{cor.converse} If $T$ is has a minimal isometric Nica-covariant dilation, then,
\begin{equation*}
\sum_{\substack{U\subseteq W \\ U\mbox{ is a clique}}} (-1)^{|U|} T_U T_U^*\geq 0.
\end{equation*} 
\end{corollary} 

\begin{proof} Let $V:P_\Gamma\to\bh{K}$ be the minimal Naimark dilation for $T$. We have $\mathcal{H}$ is co-invariant for $V$, and thus with respect to the decomposition $\mathcal{K}=\mathcal{H}\oplus\mathcal{H}^\perp$, $V(p)=\begin{bmatrix} T(p) & 0 \\ * & * \end{bmatrix}$. Therefore, for every clique $U$ in $\Gamma$, $$T_U T_U^* = P_\mathcal{H} V(e_U) V(e_U)^* \big|_\mathcal{H}.$$ 

It suffices to show for every $W\subseteq\Lambda$,
\begin{equation}\label{eq.converse}
\sum_{\substack{U\subseteq W \\ U\mbox{ is a clique}}} (-1)^{|U|} V(e_U) V(e_U)^*\geq 0.
\end{equation} 

For each vertex $i\in\Lambda$, denote $P_i=V(e_i)V(e_i)^*$ the range projection of the isometry $V(e_i)$. Since $V$ is Nica-covariant, $P_i,P_j$ commutes and 
$$P_i P_j=\begin{cases} V_i V_j V_j^* V_i^*, \mbox{ if }i\mbox{ is adjacent to }j; \\ 0, \mbox{ otherwise.}\end{cases}$$

For each $U\subseteq W$, denote $P_U=\prod_{i\in U} P_i$ and in particular let $P_\emptyset = I$. If $U\subseteq W$ is not a clique, then we can find two vertices $i,j\in U$ that are not adjacent. Since $P_iP_j=0$, it follows that $P_U=0$. If $U\subseteq W$ is a clique, then it follows from that Nica-covariant condition that $P_U=V(e_U) V(e_U)^*$. 

Consider the projection $R=\prod_{i\in W} (I-P_i)$:
\begin{align*}
R &= \prod_{i\in W} (I-P_i) \\
&= \sum_{U\subseteq W} (-1)^{|U|} P_U \\
&= \sum_{\substack{U\subseteq W \\ U\mbox{ is a clique}}} (-1)^{|U|} P_U \\
&= \sum_{\substack{U\subseteq W \\ U\mbox{ is a clique}}} (-1)^{|U|} V(e_U) V(e_U)^*.
\end{align*}

Since $R$ is a projection, $R\geq 0$ and this proves the Condition \pref{eq.converse}. \end{proof} 

We have now established the equivalence among Condition \pref{eq.main1}, $\ast$-regular, and having a minimal isometric Nica-covariant dilation.

\begin{thm:main1} Let $T:P_\Gamma\to\bh{H}$ be a representation. Then the following are equivalent:
\begin{enumerate}
\item\label{thm.main1.1} $T$ is $\ast$-regular,
\item\label{thm.main1.2} $T$ has a minimal isometric Nica-covariant dilation,
\item\label{thm.main1.3} $T$ satisfies Condition \pref{eq.main1}.
\end{enumerate}
\end{thm:main1} 

\begin{proof} \pref{thm.main1.1} $\Longrightarrow$ \pref{thm.main1.2} is established in Theorem \ref{thm.nc}. \pref{thm.main1.1} $\Longrightarrow$ \pref{thm.main1.2} is established in Corollary \ref{cor.converse}. Finally, \pref{thm.main1.3} $\Longrightarrow$ \pref{thm.main1.1} is established in Theorem \ref{thm.main}. 
\end{proof} 

\section{The Property (P)}

Popescu \cite{Popescu1999} first studied the noncommutative Poisson transform associated to a certain class of operators that satisfies the property (P). The property (P) has recently been generalized to higher rank graphs \cite{Skalski2009, Skalski2010}. It turns out that the class of operators Popescu studied can be viewed as a representation of a graph product of $\mathbb{N}$, and we thereby extend the Property (P) to representations of graph products of $\mathbb{N}$. This section proves that $\ast$-regular condition implies the property (P), and they are equivalent under certain conditions. 

Throughout this section, we fix a finite simple graph $\Gamma$ whose vertex set is denoted by $\Lambda$. 

\begin{definition}\label{df.propertyP} A contractive representation $T:P_\Gamma\to\bh{H}$ is said to have \emph{the Property (P)} if there exists $0\leq \rho<1$ so that for all $\rho\leq r\leq 1$, 

\begin{equation}\label{eq.propertyP}
\sum_{\substack{U\subseteq \Lambda \\ U\mbox{ is a clique}}} (-1)^{|U|} r^{|U|} T(e_U) T(e_U)^*\geq 0.
\end{equation} 
\end{definition} 

\begin{example} Let $\Gamma$ be a complete $k$-partite graph $K_{n_1,n_2,\cdots,n_k}$. In other words, denote $\Lambda=\{(i,j):1\leq i\leq k, 1\leq j\leq n_i\}$ be the vertex set, and $(i_1,j_1)$ is adjacent to $(i_2,j_2)$ in $\Gamma$ if and only if $i_1\neq i_2$. A contractive representation $T$ of this graph semigroup $P_\Gamma$ is uniquely determined by $T_{i,j}=T(e_{i,j})$. Here, for each $i$, $T_{i,1},\cdots,T_{i,n_i}$ are not necessarily commuting contractions. However, for each $i_1\neq i_2$, $T_{i_1,j_1}$ commutes with $T_{i_2,j_2}$. 

In \cite{Popescu1999}, Popescu considered such class of operators $\{T_{i,j}\}$ where for each $i$, $\{T_{i,j}\}_{j=1}^{n_j}$ forms a row contraction in the sense that, $$\sum_{j=1}^{n_i} T_{i,j} T_{i,j}^* \leq I.$$

This family of operators is also considered in many subsequent papers on non-commutative polyballs (see also \cite{Popescu2015, Popescu2016b}). For such family of operators, Popescu says it has the property (P) if the Condition \pref{eq.propertyP} is satisfied. It is observed in \cite{Popescu1999} that the property (P) allows one to obtain a Poisson transform and subsequently a dilation of the family of operators $\{T_{i,j}\}$.

One may observe that Definition \ref{df.propertyP} of the property (P) does not require the row contractive condition. Instead, this paper mostly considers a contractive representation $T$ of the graph product $P_\Gamma$ that satisfies the Condition \pref{eq.main1} and thus has a $\ast$-regular dilation. The row contractive condition is embedded in the Condition \pref{eq.main1}.
\end{example}

Our first result shows that if $T$ satisfies the Condition \pref{eq.main1}, then it has the property (P). Let $T:P_\Gamma\to\bh{H}$ be a representation that satisfies the Condition \pref{eq.main1}. By Theorem \ref{thm.main.all}, it has a minimal isometric Nica-covariant dilation $V:P_\Gamma\to\bh{K}$. Moreover, $\mathcal{H}$ is co-invariant for $V$, and thus $$P_\mathcal{H} V(e_U) V(e_U)^* \big|_\mathcal{H} = T(e_U)T(e_U)^*.$$

Therefore, to show $T$ has the property (P), it suffices to show $V$ has the property (P). For $r\in\mathbb{R}$, let us denote $$f(r)=\sum_{\substack{U\subseteq \Lambda \\ U\mbox{ is a clique}}} (-1)^{|U|} r^{|U|} V(e_U) V(e_U)^*.$$

It follows from the proof of Corollary \ref{cor.converse} that $f(1)\geq 0$. In fact, $f(1)$ is a projection onto the subspace that is orthogonal to all the ranges of $V(e_i)$. Following the notation we used in the proof of Corollary \ref{cor.converse}, for each vertex $i\in\Lambda$, denote $P_i=V_i V_i^*$. Since $V$ is Nica-covariant, $P_i, P_j$ commute, and $$P_i P_j=\begin{cases} V_i V_j V_j^* V_i^*, \mbox{ if }i\mbox{ is adjacent to }j; \\ 0, \mbox{ otherwise.}\end{cases}$$

For each $U\subseteq \Lambda$, denote $P_U=\prod_{i\in U} P_i$, the projection onto the intersection of the ranges of all $\{P_i\}_{i\in U}$. In particular, we let $P_\emptyset=I$. Notice that if there are two vertices $i,j\in U$ that are not adjacent, $P_i P_j=0$ and thus $P_U=0$. Therefore, $P_U\neq 0$ only if $U$ is a clique. The function $f(r)$ can be rewritten as 
\begin{align*}
f(r) &=\sum_{\substack{U\subseteq \Lambda \\ U\mbox{ is a clique}}} (-1)^{|U|} r^{|U|} P_U \\
&= \sum_{U\subseteq \Lambda} (-1)^{|U|} r^{|U|} P_U \\
&= \sum_{k=0}^{|\Lambda|} \left(\sum_{\substack{U\subseteq \Lambda \\ |U|=k}} (-1)^k  P_U \right) r^k
\end{align*}

For each $U\subseteq\Lambda$, denote $R_U=P_U\cdot \prod_{i\notin U} P_i^\perp$. The range of $R_U$ are those vectors that are contained in the range of $P_U$ but orthogonal to the range of $P_i$ where $i\notin U$. In particular, $R_\emptyset=\prod_{i\in\Lambda} P_i^\perp$, which is the projection onto those vectors that are orthogonal to the ranges of all $P_i$. It was observed in Corollary \ref{cor.converse} that $$R_\emptyset=\sum_{\substack{U\subseteq \Lambda \\ U\mbox{ is a clique}}} (-1)^{|U|} V(e_U) V(e_U)^*=f(1).$$

Finally, denote 
\begin{equation}\label{eq.Qm}
Q_m = \sum_{\substack{U\subseteq \Lambda \\ |U|=m}} R_U.
\end{equation}

In particular, $Q_0=R_\emptyset=f(1)$. Notice that if two distinct subsets $U_1,U_2\subseteq\Lambda$ and $|U_1|=|U_2|=m$, then at least one vertex in $U_1$ is not in $U_2$ and vice versa. Therefore, $R_{U_1} R_{U_2}=0$ and thus $R_{U_1},R_{U_2}$ are projections onto orthogonal subspaces. Hence, $Q_m$ is a projection. Intuitively, the range of $Q_m$ are those vectors that are contained in the range of $m$ of $P_i$ and orthogonal to the range of all other $P_i$. Therefore, $\{Q_m\}_{m=0}^{|\Lambda|}$ are pairwise orthogonal projections and $$\sum_{m=0}^{|\Lambda|} Q_m=I.$$

We first obtain a Taylor expansion of $f$ about $r=1$. For each $1\leq m\leq |\Lambda|$, the $m$-th derivative of $f$ is equal to:
\begin{align*}
f^{(m)}(r) &= \sum_{k=m}^{|\Lambda|} \sum_{\substack{U\subseteq \Lambda \\ |U|=k}} (-1)^k \frac{k!}{(k-m)!}r^{k-m} P_U \\
&= (-1)^m m! \sum_{k=m}^{|\Lambda|} \sum_{\substack{U\subseteq \Lambda \\ |U|=k}} (-1)^{k-m} {k\choose m} r^{k-m} P_U \\
\end{align*} 

\begin{lemma}\label{lm.taylorP} $f^{(m)}(1)=(-1)^m m!\cdot Q_m$. Moreover, $f$ has the Taylor series expansion $$f(r)= \sum_{m=0}^{|\Lambda|} (-1)^m (r-1)^m Q_m.$$ 
\end{lemma} 

\begin{proof} It suffices to prove $$Q_m=\sum_{k=m}^{|\Lambda|} \sum_{\substack{U\subseteq \Lambda \\ |U|=k}} (-1)^{k-m} {k\choose m} P_U.$$

Denote the right hand side of the summation $S_m$. It suffices to prove 
$$S_m Q_i=Q_i S_m=\begin{cases} 
Q_i, \mbox{ if }i=m; \\
0, \mbox{ if }i\neq m.
\end{cases}$$

From Equation \pref{eq.Qm}, $Q_m$ is the sum of all $R_W$ where $|W|=m$. Since $\{R_W\}_{|W|=m}$ are pairwise orthogonal projections, it suffices to prove
$$S_m R_W=R_W S_m=\begin{cases} 
R_W, \mbox{ if }|W|=m; \\
0, \mbox{ if }|W|\neq m.
\end{cases}$$

First of all, since $\{P_i\}_{i\in\Lambda}$ are commuting orthogonal projections, $R_W, S_m$ commute for all $W\subseteq\Lambda$ and $0\leq m\leq |\Lambda|$. Fix $W$ and consider $S_m R_W$. 

If $|W|<m$, then every $|U|\geq m$ contains some vertex not in $W$. Therefore, $P_U R_W=0$, and hence $S_m R_W=0$. 

If $|W|\geq m$, then for each $|U|\geq m$, $$P_U R_W=\begin{cases} R_W,\mbox{ if } U\subseteq W; \\ 0,\mbox{ otherwise}.\end{cases}$$

Therefore,
\begin{align*}
S_m R_W &= \left(\sum_{k=m}^{|\Lambda|} \sum_{\substack{U\subseteq \Lambda \\ |U|=k}} (-1)^{k-m} {k\choose m} P_U\right) \cdot R_W \\
&= \sum_{k=m}^{|W|} \sum_{\substack{U\subseteq W \\ |U|=k}} (-1)^{k-m} {k\choose m} R_W \\
&= \sum_{k=m}^{|W|} (-1)^{k-m} {|W| \choose k} {k\choose m}  {k\choose m} R_W \\
&= \sum_{k=m}^{|W|} (-1)^{k-m} \frac{|W|!}{k!(|W|-k)!} \frac{k!}{m!(k-m)!} R_W \\
&= {|W| \choose m} \sum_{k=m}^{|W|} (-1)^{k-m} {{|W|-m} \choose {k-m}}  R_W \\
&= {|W| \choose m} \sum_{j=0}^{|W|-m} (-1)^j {{|W|-m} \choose j}  R_W.
\end{align*}

Here, $\sum_{j=0}^{|W|-m} (-1)^j {{|W|-m} \choose j}$ is equal to $(1-1)^{|W|-m}=0$ if $|W|>m$, and $1$ if $|W|=m$. Therefore, $$S_m R_W=\begin{cases} R_W, \mbox{ if } |W|=m; \\ 0, \mbox{ otherwise}.\end{cases}$$

This proves $S_m=Q_m$. Since the graph $\Gamma$ is assumed to be a finite graph, $f(r)$ is a finite operator-valued polynomial. Its Taylor series expansion about $1$ is equal to:
\begin{align*}
f(r) &= \sum_{m=0}^{|\Lambda|} \frac{f^{(m)}(1)}{m!} (r-1)^m \\
&= \sum_{m=0}^{|\Lambda|} (-1)^m (r-1)^m Q_m. \qedhere
\end{align*}
 \end{proof} 

\begin{theorem}\label{thm.propertyP} If a representation $T:P_\Gamma\to\bh{H}$ is $\ast$-regular, then $T$ satisfies the property (P). Moreover, the constant $\rho$ in the property (P) can be chosen to be $\rho=0$.
\end{theorem} 

\begin{proof} Let $V:P_\Gamma\to\bh{K}$ be the minimal isometric $\ast$-regular dilation for $T$. By Lemma \ref{lm.taylorP}, for each $0\leq r\leq 1$,
\begin{align*}
f(r) &= \sum_{\substack{U\subseteq \Lambda \\ U\mbox{ is a clique}}} (-1)^{|U|} r^{|U|} P_U \\
&= \sum_{m=0}^{|\Lambda|} (-1)^m (r-1)^m Q_m
\end{align*}

For $0\leq r\leq 1$, $(-1)^m (r-1)^m\geq 0$. Since each $Q_m$ is an orthogonal projection, $f(r)\geq 0$. Notice when $U$ is a clique, $P_U=V_UV_U^*$, where $V_U=\begin{bmatrix} T_U & 0 \\ * & * \end{bmatrix}$ with respect to $\mathcal{K}=\mathcal{H}\oplus\mathcal{H}^\perp$. Therefore, by projecting onto the corner corresponding to $\mathcal{H}$, we obtain that for all $0\leq r\leq 1$, $$\sum_{\substack{U\subseteq \Lambda \\ U\mbox{ is a clique}}} (-1)^{|U|} r^{|U|} T_U T_U^*\geq 0.$$ 

This implies $T$ satisfies the property (P) with $\rho=0$. \end{proof} 

It is not clear when the converse of Theorem \ref{thm.propertyP} also holds. Popescu established in \cite[Corollary 5.2]{Popescu1999} the converse for a special class of operators.  

\begin{proposition}[Corollary 5.2, \cite{Popescu1999}]\label{prop.Popescu} Let $\Gamma=K_{n_1,\cdots,n_k}$ be a complete $k$-multipartite graph. Let $\{T_{i,j}\in\bh{H}: 1\leq i\leq k, 1\leq j\leq n_i\}$ be a family of operators such that:
\begin{enumerate}
\item\label{cond.P1} For each $i$, $\sum_{j=1}^{n_i} T_{i,j} T_{i,j}^*\leq I$,
\item\label{cond.P2} The associated representation $T:P_\Gamma\to\bh{H}$ has the property (P).
\end{enumerate}

Then the associated representation $T$ has a minimal isometric Nica-covariant dilation. 
\end{proposition} 

However, for a representation of an arbitrary graph semigroup, it is not clear how one can replace the Condition \pref{cond.P1} in the Proposition \ref{prop.Popescu}. 

\begin{example} Let us consider the special case when $n_1=\cdots=n_k=1$ and the graph $\Gamma$ is the complete graph on $k$-vertices. Let $\{T_i\}_{i=1}^k$ be a family of operators as in the Proposition \ref{prop.Popescu}. Notice the Condition \pref{cond.P1} is simply saying that each $T_i$ is a contraction. Proposition \ref{prop.Popescu} states that such $T_i$ has a minimal isometric Nica-covariant dilation, and thus by the Theorem \ref{thm.main.all}, $T_i$ has to satisfy the Condition \pref{eq.main1}. Note that in a complete graph, Condition \ref{eq.main1} is the same as the Brehmer's Condition \pref{eq.Brehmer.regular}.

In fact, we can derive the Condition \pref{eq.main1} directly from the property (P), without invoking the minimal isometric Nica-covariant dilation.

For any subset $W\subseteq\{1,2,\cdots,n\}$, denote $$\Delta_W(r)=\sum_{U\subseteq W} (-1)^{|U|} r^{|U|} T_U T_U^*.$$

The property (P) implies for some $0\leq \rho<1$ and all $\rho\leq r\leq 1$, $\Delta_{\{1,2,\cdots,n\}}(r)\geq 0$. For any $1\leq i\leq n$, let $W_i=\{1,\cdots,i-1,i+1,\cdots,n\}$. Notice that, $$\Delta_{\{1,2,\cdots,n\}}(r) = \Delta_{W_i}(r) - rT_i\Delta_{W_i}(r) T_i^*.$$

We claim that $\Delta_{W_i}(r)\geq 0$ for all $\rho\leq r< 1$. If otherwise, since $\Delta_{W_i}(r)$ is a self-adjoint operator, let $$-M=\inf\{\langle \Delta_{W_i}(r)h,h\rangle: \|h\|=1\}<0.$$

Pick a unit vector $h$ so that $-M\leq \langle \Delta_{W_i}(r)h,h\rangle< -M\cdot r$. Then,
\begin{align*}
\langle rT_i\Delta_{W_i}(r) T_i^*h,h\rangle &= r\cdot \langle \Delta_{W_i}(r) T_i^*h,T_i^*h\rangle \\
&\geq -M\cdot r.
\end{align*}

Therefore,
\begin{align*}
\langle \Delta_{\{1,2,\cdots,n\}}(r) h,h\rangle &= \langle \Delta_{W_i}(r)h,h\rangle - \langle rT_i\Delta_{W_i}(r) T_i^*h,h\rangle \\
&< -M\cdot r + M\cdot r = 0.
\end{align*}

This contradicts that $\Delta_{\{1,2,\cdots,n\}}(r)\geq 0$. Hence, we can conclude that $\Delta_{W_i}(r)\geq 0$. In other words, $\{T_1,\cdots,T_{i-1},T_{i+1},\cdots,T_n\}$ satisfies the property (P). Similarly, by removing one element each time, we obtain that for any $W\subseteq\{1,2,\cdots,n\}$, $\Delta_W(r)\geq 0$ for all $\rho\leq r<1$. In particular, let $r\to 1$, we obtain that for every $W\subseteq\{1,2,\cdots,n\}$, $$\sum_{U\subseteq W} (-1)^{|U|} T_U T_U^*\geq 0.$$

This is exactly the Condition \pref{eq.main1} on the complete graph (equivalently, the Brehmer's Condition \pref{eq.Brehmer.regular}).\end{example} 

\begin{remark} For an arbitrary graph $\Gamma$, it is not clear how we can replace the Condition \pref{cond.P2} in the Proposition \ref{prop.Popescu} to guarantee a minimal isometric Nica-covariant dilation for a representation $T:P_\Gamma\to\bh{H}$. 
\end{remark} 

\bibliographystyle{abbrv}
\bibliography{DilnGraphProd}

\end{document}